\documentclass[11pt]{article}
\usepackage{epsfig}

\def\be{\begin{equation}}
\def\ee{\end{equation}}
\def\bea{\begin{eqnarray}}
\def\eea{\end{eqnarray}}
\def\bes{\begin{eqnarray*}}
\def\ees{\end{eqnarray*}}

\def\nn{\nonumber}
\def\lb{\label}
\def\bs{\setminus}

\def\R{{\bf R}}
\def\C{{\bf C}}
\def\Z{{\bf Z}}
\def\N{{\bf N}}
\def\U{{\bf U}}

\def\Q{{\bf Q}}
\def\T{{\bf T}}
\def\CP{{\bf CP}}
\def\HP{{\bf HP}}
\def\CaP{{\bf CaP}}

\def\aa{{\alpha}}
\def\bb{{\beta}}
\def\ga{{\gamma}}

\def\ka{{\kappa}}
\def\th{{\theta}}

\def\om{{\omega}}
\def\Om{{\Omega}}
\def\ep{{\epsilon}}
\def\lm{{\lambda}}
\def\Lm{{\Lambda}}
\def\dl{{\delta}}
\def\Dl{{\Delta}}
\def\sg{{\sigma}}
\def\dm{{\diamond}}

\def\<{{\langle}}
\def\>{{\rangle}}

\def\K{{\mathcal K}}

\def\rank{{\rm rank}}

\def\Sp{{\rm Sp}}

\def\ol{\overline}

\def\hb{\vrule height0.18cm width0.14cm $\,$}

\title{The enhanced common index jump theorem for symplectic paths and non-hyperbolic closed
geodesics on Finsler manifolds}

\author{Huagui Duan$^{1}$,\thanks{Partially supported by NSFC (Nos. 11131004, 11471169),
LPMC of MOE of China and Nankai University. E-mail: duanhg@nankai.edu.cn.}
\quad Yiming Long$^{2}$,\thanks{Partially supported by NSFC (No. 11131004), MCME and LPMC of
MOE of China, Nankai University and BCMIIS of Capital Normal University. E-mail: longym@nankai.edu.cn.}
\quad Wei Wang$^{3}$ \thanks{Partially supported by NSFC (Nos. 11222105, 11431001),  E-mail: alexanderweiwang@gmail.com.}\\\\
$^{1}$ School of Mathematical Sciences and LPMC, Nankai University, Tianjin 300071\\
$^{2}$ Chern Institute of Mathematics and LPMC, Nankai University, Tianjin 300071\\
$^{3}$ School of Mathematical Sciences and LMAM, Peking University, Beijing 100871\\
The People's Republic of China\\}

\begin{document}
\date{December 31, 2015}
\maketitle

\begin{abstract}
{\it }In this paper, we first generalize the common index jump theorem for symplectic matrix paths proved
in 2002 by Long and Zhu in \cite{LoZ}, and get an enhanced version of it. As its applications, we
further prove that for a compact simply-connected manifold $(M,F)$ with a bumpy irreversible Finsler
metric $F$ and $H^*(M;\Q)\cong T_{d,n+1}(x)$ for some even integer $d\ge 2$ and integer $n\ge 1$,
there exist at least $\frac{dn(n+1)}{2}$ distinct non-hyperbolic closed geodesics with odd Morse
indices, provided the number of distinct prime closed geodesics is finite and every prime closed geodesic
satisfies $i(c)>0$. Note that the last non-zero index condition is satisfied if the flag curvature $K$
satisfies $K\ge 0$. For an odd-dimensional bumpy Finsler sphere $(S^d,F)$, there exist at least $(d+1)$
distinct prime closed geodesics with even Morse indices, and at least $(d-1)$ of which are non-hyperbolic,
provided the number of distinct prime closed geodesics is finite and every prime closed geodesic $c$
satisfies $i(c)\ge 2$. Note that the last index condition $i(c)\ge 2$ is satisfied if the reversibility
$\lm$ and the flag curvature $K$ of $(M,F)$ satisfy $\frac{\lm^2}{(1+\lm)^2}<K\le 1$. Note that the first
two in the above three lower bound estimates are sharp due to examples constructed by Katok and Ziller.
In addition, we also prove that either there exists at least one non-hyperbolic closed geodesic, or there
exist infinitely many distinct closed geodesics on a compact simply connected bumpy Finsler $(M,F)$
satisfying the above cohomological condition with some even integer $d\ge 2$ and integer $n\ge 1$.
\end{abstract}



{\bf Running head}: Closed geodesics on Finsler manifolds.

\renewcommand{\theequation}{\thesection.\arabic{equation}}
\renewcommand{\thefigure}{\thesection.\arabic{figure}}

\setcounter{figure}{0}
\setcounter{equation}{0}
\section{Introduction and main results}

The common index jump theorem (cf. Theorem 4.2 of \cite{LoZ}, which is denoted by CIJT for short below)
was established by Long and Zhu in \cite{LoZ} of 2002. Since then, it has become an important tool in
the studies of multiplicity and stability problems of closed characteristics on
given energy hypersurfaces of Hamiltonian dynamics and closed geodesics on Finsler manifolds. Denote
by $i(\ga,m)$ and $\nu(\ga,m)$ the Maslov-type index and nullity of the $m$-th iterate $\ga^m$ of any
symplectic matrix path $\ga$ as in \cite{LoZ} and \cite{Lon3}. For a finite family of symplectic paths
$\{\ga_k\}_{1\le k\le q}$ in $\Sp(2n)$ with positive mean indices $\hat{i}(\ga_k)>0$, the CIJT yields
$(q+1)$-tuple of integers $(N, m_1, \ldots, m_k)$ so that all the index jump intervals
$$ (i(\ga_k, 2m_k-1)+\nu(\ga_k,2m_k-1)-1, i(\ga_k,2m_k+1))  $$
with $1\le k\le q$ possesses a large enough common intersection interval $[2N-\ka_1, 2N+\ka_2]$ for
some positive constants $\ka_1$ and $\ka_2$. Then it was further proved that the number of integers
in $2\N-2+n$ contained in this common interval yields a lower bound on the number of distinct closed
characteristics on compact convex energy hypersurfaces in $\R^{2n}$ under study. Here one important
feature in the index jump intervals is that only indices of $(2m_k\pm 1)$-th iterates $\ga_k^{2m_k\pm 1}$
and certain estimate on the $2m_k$-th iterate $\ga_k^{2m_k}$ are needed, based on the fact that the initial
index $i(\ga_k)\ge n$ holds always, which implies the monotone increasing property of index $i(\ga_k,m)$ as
$m$ increases and the index intervals $\{[i(\ga_k,m),i(\ga_k,m)+\nu(\ga_k,m)-1]\}_{m\ge 1}$ are mutually
non-intersecting.

In this paper our main goal is to study multiplicity and stability problems of closed geodesics on
any compact simply-connected Finsler manifold $(M,F)$. For such an orbit $c$, the initial Morse index $i(c)$
can be rather small like $0$ or $1$ even if $\dim M$ is rather large, and one can not hope the monotone
increasing property of the Morse index $i(c^m)$ as $m$ increases and the index interval non-intersecting
property. Thus in order to continue to use ideas of the CIJT, we need to understand precisely the indices of
all the $(2m_k\pm m)$-th iterates $\ga_k^{2m_k\pm m}$ for every integer $1\le m \le \bar{m}$ with any
given $\bar{m}>0$ as well as the precise index of $2m_k$-th iterate $\ga_k^{2m_k}$. Thus in the first part
of this paper, we generalize the CIJT to our Theorem 3.5 below to obtain such precise values of Maslov-type
indices of these iterates of closed geodesics. Such information allows us to compute the corresponding Morse
inequality accurately and derive certain sharp estimates on the multiplicity and stability of closed
geodesics compact simply-connected Finsler manifolds as in the proof of Theorem 1.1 below. We call this
generalization of CIJT (Theorem 3.5 below) the enhanced common index jump theorem (enhanced CIJT for short
below). We believe that this enhanced CIJT can be applied to many other problems on periodic solution orbits
in Hamiltonian and symplectic dynamics too.

There is a famous conjecture in Riemannian geometry which claims the existence of infinitely
many distinct closed geodesics on every compact Riemannian manifold $M$. This conjecture has
been proved for many cases (cf. for example \cite{BTZ1}, \cite{BTZ2}, \cite{Ban1}), specially
Gromoll and Meyer proved the conjecture in \cite{GrM} of 1969 for any compact $M$ provided the
Betti number sequence $\{b_p(\Lm M)\}_{p\in\N}$ of the free
loop space $\Lm M$ of $M$ is unbounded. Then in \cite{ViS1} of 1976, for compact simply
connected manifold $M$, Vigu\'e-Poirrier and Sullivan further proved this Betti number
sequence is bounded if and only if $M$ satisfies
\be  H^*(M;\Q)\cong T_{d,n+1}(x)=\Q[x]/(x^{n+1}=0)   \lb{1.1}\ee
with a generator $x$ of degree $d\ge 2$ and height $n+1\ge 2$, where $\dim M=dn$. Note that
when $d$ is odd, then $x^2=0$ and $n=1$, or when $n=1$, $M$ is rationally homotopic to $S^d$
(cf. Remark 2.5 of \cite{Rad1} and \cite{Hin}). Among these manifolds, only for Riemannian $S^2$ this
conjecture was proved by the works \cite{Fra} of Franks in 1992 and \cite{Ban2} of Bangert in 1993.
When Finsler manifolds are considered, the situation changes dramatically. It was quite surprising
that Katok found some irreversible Finsler metrics on rank one symmetric spaces which possess
finitely many distinct closed geodesics in \cite{Kat} of 1973. There Katok carried out also
detailed constructions of the Finsler metrics on $S^d$ which possess precisely $2[(d+1)/2]$
distinct closed geodesics and all of which are non-degenerate and elliptic. The geometry of
Katok's metrics was further studied by Ziller in \cite{Zil} of 1982. He constructed in detail
Finsler metrics for complex projective spaces $\CP^n$ (with $d=2$), quaternionic
projective spaces $\HP^n$ (with $d=4$), and the Cayley plane $\CaP^2$ (with $d=8$ and $n=2$)
which possess precisely $n(n+1)$, $2n(n+1)$, and $24$ distinct closed geodesics respectively,
i.e., precisely $\frac{dn(n+1)}{2}$ distinct closed geodesics. Based on Katok's work, people
believe that the smallest number of distinct closed geodesics on any Finsler sphere
$S^d$ is $2[(d+1)/2]$, as conjectured by Anosov in \cite{Ano} of 1974. This was proved in
2005 by Bangert and Long for every Finsler sphere $(S^2,F)$ (cf. \cite{Lon4} of 2006), which
was published latter as \cite{BaL} in 2010. Now it is natural to generalize this conjecture to
all compact simply connected Finsler manifolds satisfying (\ref{1.1}), i.e., the lower bound
of the number of distinct closed geodesics on such manifolds should be $\frac{dn(n+1)}{2}$ when
$d\ge 2$ is even and $d+1$ when $d\ge 2$ is odd respectively. The second part of this paper is
devoted to study this problem as well as the non-hyperbolicity of these closed geodesics.

Recall that a closed curve on a Finsler manifold is a closed geodesic if it is locally the shortest
path connecting any two nearby points on this curve (cf. \cite{She}). As usual, on any Finsler
manifold $(M,F)$, and a closed geodesic $c:S^1=\R/\Z\to M$ is {\it prime} if it is not a multiple
covering (i.e., iteration) of any other closed geodesics. Here the $m$-th iteration $c^m$ of $c$ is
defined by $c^m(t)=c(mt)$. The inverse curve $c^{-1}$ of $c$ is defined by $c^{-1}(t)=c(1-t)$ for
$t\in \R$. Note that unlike the Riemannian case, the inverse curve $c^{-1}$ of a closed geodesic $c$
on an irreversible Finsler manifold need not be a geodesic. We call two prime closed geodesics $c$
and $d$ {\it distinct} on a Finsler manifold, if there exists no $\th\in (0,1)$ such that
$c(t)=d(t+\th)$ for all $t\in\R$. On a reversible Finsler (or Riemannian) manifold, two closed
geodesics $c$ and $d$ are called { \it geometrically distinct} if $c(S^1)\neq d(S^1)$, i.e., their
image sets in $M$ are distinct. We shall omit the word {\it distinct} when we talk about more than
one prime closed geodesic for simplicity.

For a closed geodesic $c$ on $n$-dimensional manifold $(M,\,F)$, denote by $P_c$ the linearized
Poincar\'{e} map of $c$. Then $P_c\in \Sp(2n-2)$ is well known. As usual, for any $M\in \Sp(2k)$, we
define the {\it elliptic height} $e(M)$ of $M$ to be the total algebraic multiplicity of all
eigenvalues of $M$ on the unit circle $\U=\{z\in\C|\; |z|=1\}$ in the complex plane $\C$. Since $M$
is symplectic, $e(M)$ is even and $0\le e(M)\le 2k$. A closed geodesic $c$ is {\it elliptic} if
$e(P_c)=2(n-1)$, i.e., all the eigenvalues of $P_c$ locate on $\U$; {\it irrationally elliptic} if
it is elliptic and all the eigenvalues of $P_c$ locate on $\U\bs\{\pm 1\}$; {\it hyperbolic} if
$e(P_c)=0$, i.e., all the eigenvalues of $P_c$ locate away from $\U$; {\it non-degenerate} if $1$ is
not an eigenvalue of $P_c$. A Finsler manifold $(M,\,F)$ is called {\it bumpy} if all the closed
geodesics and their iterates on $M$ are non-degenerate (cf. \cite{Abr}).

Recently the Maslov-type index theory for symplectic paths has been applied to study the closed
geodesic problem. In 2005, Bangert and Long proved the existence of at least two distinct closed
geodesics on every Finsler $(S^2,F)$ (which was published as \cite{BaL} in 2010). Since then in the
last ten years, a great number of results on the multiplicity and stability of closed geodesics on
Finsler manifolds has appeared, for which we refer readers to \cite{DuL1}-\cite{DuL3}, \cite{LoD},
\cite{LoW}, \cite{Rad4}-\cite{Rad6}, \cite{Wan1}-\cite{Wan2}, \cite{HiR}, \cite{XiL}, \cite{DLX},
and the references therein.

Recently Wang in \cite{Wan2} proved the existence of at least $2[\frac{n+1}{2}]$ distinct closed
geodesics on a bumpy Finsler $(S^n,F)$ when the flag curvature $K$ satisfies
$\frac{\lm^2}{(1+\lm)^2}<K\le 1$, where the reversibility $\lambda=\lambda(M,F)$ was introduced
by Rademacher in \cite{Rad3} as
\be   \lambda=\max\{F(-X)\ |\ X\in TM,\ F(X)=1\}\ge 1.   \lb{1.2}\ee
Note that the lower bound $2[\frac{n+1}{2}]$ on the number of distinct closed geodesics for spheres $S^n$
is sharp due to the above mentioned examples of Katok. Besides spheres, in \cite{Rad4}, Rademacher obtained some
multiplicity and stability results of closed geodesics on compact simply connected manifolds satisfying
some pinching conditions. For example, it was proved
that a Finsler $(\CP^n,F)$ with $n\ge 7$ and a bumpy Finsler metric possessing only finitely many distinct
closed geodesics and satisfying $\left(\frac{2}{n+1}\frac{\lambda}{1+\lambda}\right)^2\le K\le 1$, carries
always at least $2n$ distinct closed geodesics, and at least $(n-3)$ of them are non-hyperbolic, where
$\lambda$ and $K$ denote also the reversibility and the flag curvature of $(\CP^n, F)$ respectively. In
\cite{Rad5}, Rademacher further proved the existence of two prime closed geodesics on any $\CP^2$ with a
bumpy irreversible Finsler metric. In recent preprint \cite{DLW}, the authors proved the existence of at
least two prime closed geodesics on every compact simply-connected Finsler manifold $(M,F)$ with a bumpy
irreversible Finsler metric $F$.

Motivated by these results, in this paper we prove

\medskip

{\bf Theorem 1.1.} {\it On every compact simply-connected manifold $(M,F)$ with a bumpy, irreversible Finsler
metric $F$ satisfying $H^*(M;\Q)\cong T_{d,n+1}(x)$ for some even integer $d\ge 2$ and some integer $n\ge 1$,
there exist at least $\frac{dn(n+1)}{2}$ distinct non-hyperbolic prime closed geodesics with odd Morse indices,
provided the number of distinct prime closed geodesics is finite and every prime closed geodesics possesses
non-zero Morse index.}

\medskip

Note that the monotonicity of iterated Morse indices of closed geodesics has played an important and
crucial role in many studies on closed geodesics (cf. \cite{DuL1}, \cite{Rad6} and \cite{Wan1})-\cite{Wan2}.
For example, Rademacher in \cite{Rad6} and Wang in \cite{Wan1}-\cite{Wan2} used the monotonicity and the
common index jump theorem of Long and Zhu established in \cite{LoZ} to get some multiplicity results of
closed geodesics on Finsler spheres. However, such a monotonicity may no longer hold in general if the
initial Morse index $i(c)$ is too small, when closed geodesics on a higher dimensional compact
simply-connected manifold $(M,F)$ are studied, even if every prime closed geodesic in the study has
non-zero Morse index. Instead of such index monotonicity, the proof of our Theorem 1.1 is based on our
enhanced common index jump theorem 3.5.

For reader's convenience, we describe the main ideas of the proof of Theorem 1.1 briefly in three steps.

More precisely, in order to prove Theorem 1.1, we assume that there exist only finitely many distinct
prime closed geodesics $\{c_k\}_{k=1}^q$ with $i(c_k)\ge 1$ for $1\le k\le q$ on such an $(M,F)$.

{\bf Step 1.} Firstly, we apply the enhanced CIJT (Theorem 3.5) to get a $(q+1)$-tuple
$(N, m_1, \cdots, m_q)$ $\in\N^{q+1}$ such that the Morse indices $i(c_k^h)$ of $h=2m_k\pm m$ and $h=2m_k$-th
iterates of each prime closed geodesic $c_k$ satisfy the inequalities (\ref{4.12}), (\ref{4.14}) and the
equality (\ref{4.13}) below with $m\in [1,\bar{m}]$ for some suitably chosen $\bar{m}\in\N$. Here the
condition $i(c_k)>0$ is used for $1\le k\le q$.

Then using these information on the Morse indices, we can compute
the alternating sum of Morse type numbers up to $2N$, i.e., $\sum_{p=0}^{2N}(-1)^pM_p$, which becomes the
sum of $\sum_{k=1}^q2m_k\ga_{c_k}=2NB(d,n)$ and $(N_+^o - N_-^o)$, where
\be N_{\pm}^o = \;^{\#}\{1\le k\le q\ |\ \pm(i(c_k^{2m_k})- 2N \mp 1) \ge 0,\ i(c_k^{2m_k})-i(c_k)\in 2\N_0,\
            i(c_k)\in 2\N-1\}. \lb{1.3}\ee
Now by direct computations and the Morse inequality, we obtain
$$ 2NB(d,n) + N_+^o - N_-^o = \sum_{p=0}^{2N}(-1)^pM_p \ge \sum_{p=0}^{2N}(-1)^pb_p(\Lm M)
    = 2NB(d,n) + \frac{dn(n+1)}{4},  $$
i.e.,
\be N_+^o \ge \frac{dn(n+1)}{4}. \lb{1.4}\ee
Here $B(d,n)$, $b_p(\Lm M)$ and $\ga_{c_k}$ are defined in Section 2 below.

{\bf Step 2.} Note that the $(q+1)$-tuple $(N, m_1, \cdots, m_q)$ in Step 1 is chosen according to a vertex
$\chi$ of the cub $[0,1]^l$ given in (iii) of Remark 3.6. Then for the vertex $\hat{\chi}=\hat{1}-\chi$ opposite
to $\chi$ in $[0,1]^l$ as in Figure 4.1 below, by Theorem 3.5 again, we find another $(q+1)$-tuple
$(N', m_1', \cdots, m_q') \in\N^{q+1}$, such that the inequalities (\ref{4.26})-(\ref{4.28}) hold for any
$1\le k\le q$. By the discussion in Step 1, we obtain the corresponding numbers $N_{\pm}^{'o}$ defined in
the same way as in (\ref{1.3}), and it yields
\be  N_+^{'o} \ge \frac{dn(n+1)}{4}. \lb{1.5}\ee
Then the symmetry of $\chi$ and $\hat{\chi}$ yields $N_-^o = N_+^{'o}$. Thus we obtain
$$  q \ge N_+^o + N_-^o \ge \frac{dn(n+1)}{4} + \frac{dn(n+1)}{4} = \frac{dn(n+1)}{2}.  $$
That is, the total number of distinct closed geodesics on such manifold in Theorem 1.1 is at least
$\frac{dn(n+1)}{2}$.

{\bf Step 3.} Then the non-hyperbolicity of these closed geodesics can be easily obtained according to their
index information given in the definition (\ref{1.3}) of $N_+^o$ and $N_-^o$, because if $c_k$ is hyperbolic,
it must satisfy $i(c_k^{2m_k})$ being even.

\medskip

In order to understand the non-zero Morse index condition in Theorem 1.1, let $(M,F)$ be an even
dimensional compact simply-connected Finsler manifold and $c$ be a prime closed geodesic on it. Denote
by
\be  \Pi_{c}: T_{c(0)}M\rightarrow T_{c(0)}M  \lb{1.6}\ee
the parallel transport along the geodesic $c(t)$ for $0\le t\le 1$ and $T^{\perp}$ the
$g_T$-orthogonal complement of the unit vector $T=\dot{c}(0)$ in $T_{c(0)}M$, where we use notations
in the proof of Theorem 8.8.1 of \cite{BCS}. Since $\Pi_c$ preserves the $g_T$ lengths and $g_T$
angles, $\Pi_c: T^{\perp}\rightarrow T^{\perp}$ is well-defined and is an orthogonal transformation
with $\det(\Pi_c)=1$. Since $\dim M \in 2\N$, $\dim T^{\perp}$ is odd. It then yields that
$1\in\sigma(\Pi_c)$. Therefore there is a vector $u_c\in T_{c(0)}M$ which is $g_T$-orthogonal to
$T=\dot{c}(0)$ and $\Pi_c(u_c)=u_c$. For any eigenvector $u$ belonging to the eigenvalue $1$ of
$\Pi_c$, denote by $U_{\dot{c},u}(t)$ the parallel transport of $u$ along $c(t)$ for $0\le t\le 1$
with $g_{\dot{c}}(U_{\dot{c},u},U_{\dot{c},u})=1$. As in Section 2 of \cite{Rad3}, denote the flag
curvature of $(M,F)$ at the plane $V_{\dot{c},U_{\dot{c},u}}$ spanned by $\dot{c}$ and $U_{\dot{c},u}$
by $K(\dot{c},U_{\dot{c},u})=g_{\dot{c}}(R^{\dot{c}}(\dot{c},U_{\dot{c},u})U_{\dot{c},u},\dot{c})$.
Then we introduce

{\bf Definition 1.2.} {\it The minimal and maximal average flag curvatures $\bar{K}_-(M,F)$ and
$\bar{K}_+(M,F)$ of $(M,F)$ are defined respectively by
\bea
\bar{K}_-(M,F) &=& \inf\{\int_0^1K(\dot{c},U_{\dot{c},u})(t)dt\;|\;c\;{\it is\;a\;prime\;closed\;geodesic\;on}\;M,  \nn\\
               & & \qquad\qquad\;u\;{\it is\;an\;eigenvector\;belonging\;to\;the\;eigenvalue}\;1\;{\it of}\;\Pi_c\},  \nn\\
\bar{K}_+(M,F) &=& \sup\{\int_0^1K(\dot{c},U_{\dot{c},u})(t)dt\;|\;c\;{\it is\;a\;prime\;closed\;geodesic\;on}\;M,  \nn\\
               & & \qquad\qquad\;u\;{\it is\;an\;eigenvector\;belonging\;to\;the\;eigenvalue}\;1\;{\it of}\;\Pi_c\}.  \nn\eea
If there exists no closed geodesics on $(M,F)$, we set $\bar{K}_-(M,F) = \bar{K}_+(M,F) = +\infty$.}

By our above discussions, $\bar{K}_-(M,F)$ and $\bar{K}_+(M,F)$ are well-defined when $M$ is compact, simply-connected
and $\dim M \in 2\N$. Note that specially $K\ge 0$ implies $\bar{K}_-(M,F)\ge 0$.

Note that our Theorem 1.3 below realizes the non-zero index condition in Theorem 1.1.

\medskip

{\bf Theorem 1.3.} {\it Let $(M,F)$ be an even dimensional compact simply-connected manifold with a
bumpy irreversible Finsler metric $F$. Then every prime closed geodesic $c$ on $(M,F)$ satisfies
$i(c)>0$, provided $\bar{K}_-(M,F)\ge 0$. This last condition holds specially when the flag curvature
$K\ge 0$.

Consequently, for any compact simply-connected manifold $(M,F)$ with a bumpy, irreversible Finsler
metric $F$ and $H^*(M;\Q)\cong T_{d,n+1}(x)$ for some even integer $d\ge 2$ and some integer $n\ge 1$,
if the flag curvature $K\ge 0$, then either there exist infinitely many distinct prime closed geodesics,
or there exist at least $\frac{dn(n+1)}{2}$ distinct non-hyperbolic prime closed geodesics with odd
Morse indices.}

\medskip

{\bf Remark 1.4.} Note that the manifold $M$ in Theorem 1.1 and Theorem 1.3 includes the even-dimensional
spheres $S^d$ with $n=1$ as a special case. Thus our results improved the main Theorem 1.2 for the
even-dimensional spheres $S^d$ in \cite{Wan2} and the main Theorem 1.2 of \cite{Wan3} where the stronger
index restriction $i(c)\ge d-1$ or the classical curvature pinching condition
$\left(\frac{\lm}{1+\lm}\right)^2<K\le 1$ is assumed.

\medskip

Next we study the multiplicity and stability of closed geodesics on odd-dimensional spheres $S^d$
with $n=1$ in (\ref{1.1}), which is not included in manifolds considered by Theorems 1.1 and 1.3.

Here for the odd-dimensional sphere $S^d$ our two theorems below improve the main Theorem 1.2 of
\cite{Wan2}, where the stronger index restriction $i(c)\ge d-1$ or the curvature pinching condition
$\left(\frac{\lm}{1+\lm}\right)^2<K\le 1$ is used.

\medskip

{\bf Theorem 1.5.} {\it Let $(S^d,F)$ be a bumpy Finsler sphere with an odd integer $d\ge 3$. If the
number of prime closed geodesics is finite and every prime closed geodesic $c$ satisfies $i(c)\ge 2$,
then there exist at least $(d+1)$ prime closed geodesics with even Morse indices, and $(d-1)$ ones of
them are non-hyperbolic.}

\medskip

{\bf Theorem 1.6.} {\it Let $(S^d,F)$ be a bumpy Finsler sphere with an odd integer $d\ge 3$. If
the number of prime closed geodesics is finite and the flag curvature $K$ satisfies
$\frac{\lm^2}{(1+\lm)^2}<K\le 1$, then every prime closed geodesic $c$ satisfies $i(c)\ge 2$.
Consequently by Theorem 1.5 there exist at least $(d+1)$ prime closed geodesics with even
Morse indices, and $(d-1)$ ones of them are non-hyperbolic.}

\medskip

{\bf Remark 1.7.} (i) For an odd-dimensional Finsler sphere $(S^{d},F)$, the lower bound of curvature in
the pinching condition in Theorem 1.6 can not be relaxed to smaller than $\frac{\lm^2}{(1+\lm)^2}$
due to some geometric reasons. And the multiplicity result in Theorem 1.6 under this curvature
pinching condition was used to get Theorem 1.2 of \cite{Wan2}. But the parity of the Morse indices
and the non-hyperbolicity of these closed geodesics are new.

(ii) We would like to draw readers' attentions to the remarkable paper \cite{HiR}, in which among other
results Hingston and Rademacher proved the existence of at least two distinct closed geodesics on a
sphere $(S^n,F)$ when the flag curvature satisfies $\left(\frac{\lm}{1+\lm}\right)^2 < K \le 1$.

\medskip

If the index restrictions and the curvature pinching conditions are given up, the following result
can still hold.

\medskip

{\bf Theorem 1.8.}  {\it On every compact simply-connected manifold $(M,F)$ with a bumpy, irreversible
Finsler metric $F$ satisfying (\ref{1.1}) for some even integer $d\ge 2$ and some integer $n\ge 1$,
either there exist at least one non-hyperbolic closed geodesic, or there exist infinitely many distinct
closed geodesics.}

\medskip

{\bf Remark 1.9.} (i) By the Finsler metrics constructed by Katok in \cite{Kat} and Ziller in \cite{Zil},
the lower bounds on the numbers of distinct prime closed geodesics in Theorems 1.1, 1.3, 1.5, and 1.6
are sharp. Note that bumpy is a generic condition, and our above theorems give generic results on the
above mentioned manifolds.

(ii) Based on our above theorems, the works \cite{HWZ1} in 1998 and \cite{HWZ2} in 2003 of Hofer, Wysocki
and Zehnder, \cite{BaL} in 2010 of Bangert and Long, \cite{LoW} in 2007 of Long and Wang, and other
results mentioned above, we believe that it is natural to propose the following

{\bf Conjecture.} {\it On every compact simply connected Finsler manifold $(M,F)$ satisfying (\ref{1.1}),
there exist either infinitely many distinct prime closed geodesics, or there exist precisely
$\frac{dn(n+1)}{2}$ (when $d\ge 2$ is even) or $(d+1)$ (when $d\ge 2$ is odd) distinct prime closed
geodesics. Moreover, when the total number of prime closed geodesics is finite, all of them must
be irrationally elliptic.}

\medskip

In this paper, we denote by $\N$, $\N_0$, $\Z$, $\Q$, $\R$, and $\C$ the sets of natural integers,
non-negative integers, integers, rational numbers, real numbers, and complex numbers respectively.
We use only singular homology modules with $\Q$-coefficients. For any $a\in\R$, we use the functions
\be  [a]=\max\{k\in\Z\,|\,k\le a\}, \quad
           E(a)=\min\{k\in\Z\,|\,k\ge a\}, \quad \{a\}=a-[a].  \lb{1.7}\ee

\setcounter{figure}{0}
\setcounter{equation}{0}
\section{Morse theory of closed geodesics}

Let $M=(M,F)$ be a compact Finsler manifold, the space $\Lambda=\Lambda M$ of $H^1$-maps
$\gamma:S^1\rightarrow M$ has a natural structure of Riemannian Hilbert manifolds on which the
group $S^1=\R/\Z$ acts continuously by isometries (cf. \cite{Kli}). This action is defined by
$(s\cdot\ga)(t)=\ga(t+s)$ for all $\gamma\in\Lm$ and $s, t\in S^1$. For any $\gamma\in\Lambda$,
the energy functional is defined by
\be E(\gamma)=\frac{1}{2}\int_{S^1}F(\gamma(t),\dot{\gamma}(t))^2dt.  \lb{2.1}\ee
It is $C^{1,1}$ and invariant under the $S^1$-action. The critical points of $E$ of positive
energies are precisely the closed geodesics $\gamma:S^1\to M$. The index form of the functional
$E$ is well defined on any closed geodesic $c$ on $M$, which we denote by $E''(c)$. As usual,
we denote by $i(c)$ and $\nu(c)$ the Morse index and nullity of $E$ at $c$ (cf. \cite{Kli} and
\cite{Mor1}). In the following, we denote by
\be \Lm^\kappa=\{d\in \Lm\;|\;E(d)\le\kappa\},\quad \Lm^{\kappa-}=\{d\in \Lm\;|\; E(d)<\kappa\},
  \quad \forall \kappa\ge 0. \nn\ee
For a closed geodesic $c$ we set $ \Lm(c)=\{\ga\in\Lm\mid E(\ga)<E(c)\}$.

For $m\in\N$ we denote the $m$-fold iteration map
$\phi_m:\Lambda\rightarrow\Lambda$ by $\phi_m(\ga)(t)=\ga(mt)$, for all
$\,\ga\in\Lm, t\in S^1$, as well as $\ga^m=\phi_m(\gamma)$. If $\gamma\in\Lambda$
is not constant then the multiplicity $m(\gamma)$ of $\gamma$ is the order of the
isotropy group $\{s\in S^1\mid s\cdot\gamma=\gamma\}$. For a closed geodesic $c$,
the mean index $\hat{i}(c)$ is defined as usual by
$\hat{i}(c)=\lim_{m\to\infty}i(c^m)/m$. Using singular homology with rational
coefficients we consider the following critical $\Q$-module of a closed geodesic
$c\in\Lambda$:
\be \overline{C}_*(E,c)
   = H_*\left((\Lm(c)\cup S^1\cdot c)/S^1,\Lm(c)/S^1\right). \lb{2.3}\ee

{\bf Lemma 2.1.} (cf. Satz 6.11 of \cite{Rad2} ) {\it Let $c$ be a
prime closed geodesic on a bumpy Finsler manifold $(M,F)$. Then there holds}
$$ \overline{C}_q( E,c^m) = \left\{\begin{array}{cc}
     \Q, &\quad {\it if}\;\; i(c^m)-i(c)\in 2\Z\;\;{\it and}\;\;
                   q=i(c^m),\;  \cr
     0, &\quad {\it otherwise}. \end{array}\right.  $$

{\bf Definition 2.2.} (cf. Definition 1.6 of \cite{Rad1}) {\it For a
closed geodesic $c$, let $\ga_c\in\{\pm\frac{1}{2},\pm1\}$ be the
invariant defined by $\ga_c>0$ if and only if $i(c)$ is even, and
$|\ga_c|=1$ if and only if $i(c^2)-i(c)$ is even. }

{\bf Theorem 2.3.} (cf. Theorem 3.1 of \cite{Rad1} and Satz 7.9 of \cite{Rad2}) {\it Let
$(M,F)$ be a compact simply connected bumpy Finsler manifold with
$\,H^{\ast}(M,\Q)=T_{d,n+1}(x)$. Denote prime closed geodesics on $(M,F)$
with positive mean indices by $\{c_j\}_{1\le j\le q}$ for some $q\in\N$.
Then the following identity holds
\be \sum_{j=1}^q\frac{\ga_{c_j}}{\hat{i}(c_j)}=B(d,n)
=\left\{\begin{array}{cc}
     -\frac{n(n+1)d}{2d(n+1)-4}, &\quad d\;\;{\it is\;even},\cr
     \frac{d+1}{2d-2}, &\quad d\;\;{\it is\;odd}, \end{array}\right.   \lb{2.4}\ee
where $\dim M=dn$. Here $n=1$ holds when $M$ is a sphere $S^d$ of dimension $d$.}

\medskip

{\bf Lemma 2.4.} (cf. Theorem 2.4 of \cite{Rad1} and Lemma 2.5 and Lemma 2.6 of \cite{DuL2})
{\it Let $M$ be a compact simply connected manifold with $H^*(M;\Q)\cong T_{d,n+1}(x)$ for
some integers $d\ge 2$ and $n\ge 1$.

(i) When $d$ is odd (which implies that $n=1$), i.e. M is rationally homotopic to the odd
dimensional sphere $S^d$, then the Betti numbers of the free loop space of $S^d$ are given
by
\bea b_p(\Lm M)
&=& \rank H_p(\Lm S^d/S^1,\Lm^0 S^d/S^1;\Q)  \nn\\
&=& \left\{\begin{array}{ccc}
    2,&\quad {\it if}\quad p\in \K\equiv \{j(d-1)\,|\,2\le j\in\N\},  \cr
    1,&\quad {\it if}\quad p\in \{d-1+2j\,|\,j\in\N_0\}\bs\K,  \cr
    0 &\quad {\it otherwise}. \end{array}\right. \lb{2.5}\eea
For any integer $l\ge d-1$, there holds
\bea \sum_{p=0}^lb_p(\Lm M) = \left[\frac{l}{d-1}\right] + \left[\frac{l}{2}\right] - \frac{d-1}{2}.   \lb{2.6}\eea

 (ii) When $d$ is even, let $D=d(n+1)-2$ and
\bea \Om(d,n) = \{k\in 2\N-1&\,|\,& iD\le k-(d-1)=iD+jd\le iD+(n-1)d\;  \nn\\
         && \mbox{for some}\;i\in\N\;\mbox{and}\;j\in [1,n-1]\}. \nn\eea
Then the Betti numbers of the free loop space of $M$ are given by
\bea b_p(\Lm M)
&=& \rank H_p(\Lm M/S^1,\Lm^0 M/S^1;\Q)\nn\\
&=& \left\{\begin{array}{cccc}
    0, & \quad \mbox{if}\ p\ \mbox{is even or}\ p\le d-2,  \cr
    \left[\frac{p-(d-1)}{d}\right]+1, & \quad \mbox{if}\ p\in 2\N-1\;\mbox{and}\;d-1\le p < d-1+(n-1)d, \cr
    n+1, & \quad \mbox{if}\ p\in \Om(d,n), \cr
    n, & \quad \mbox{otherwise}. \end{array}\right.\lb{2.7}\eea
For any integer $l\ge dn-1$, we have
\be \sum_{p=0}^lb_p(\Lm M) = \frac{n(n+1)d}{2D}(l-(d-1)) - \frac{n(n-1)d}{4} + 1 + \ep_{d,n}(l),  \lb{2.8}\ee
where
\bea \ep_{d,n}(l)
&=& \left\{\frac{D}{dn}\left\{\frac{l-(d-1)}{D}\right\}\right\}
          - \left(\frac{2}{d}+\frac{d-2}{dn}\right)\left\{\frac{l-(d-1)}{D}\right\}   \nn\\
&&\quad - n\left\{\frac{D}{2}\left\{\frac{l-(d-1)}{D}\right\}\right\}
          - \left\{\frac{D}{d}\left\{\frac{l-(d-1)}{D}\right\}\right\}. \lb{2.9}\eea}

\medskip

{\bf Theorem 2.5.} (cf. Theorem I.4.3 of \cite{Cha}) {\it Suppose that there exist only finitely many
prime closed geodesics $\{c_k\}_{1\le k\le q}$ on a Finsler manifold $(M,F)$. Define
\be M_p=\sum_{1\le k\le q,\ m\ge 1}\dim \ol{C}_p(E,c_k^m),\qquad\forall\ p\in\Z. \lb{2.10}\ee
Then the following Morse inequality holds for every $p\in \N_0$,
\bea M_p &\ge& b_p(\Lm M),\lb{2.11}\\
M_p - M_{p-1} + \cdots +(-1)^{p}M_0
&\ge& b_p(\Lm M) - b_{p-1}(\Lm M)+ \cdots + (-1)^{p}b_0(\Lm M).  \lb{2.12}\eea}

\setcounter{figure}{0}
\setcounter{equation}{0}
\section{The enhanced common index jump theorem of symplectic paths}

In \cite{Lon1} of 1999, Y. Long established the basic normal form decomposition of symplectic matrices.
Based on this result he further established the precise iteration formulae of indices of symplectic
paths in \cite{Lon2} of 2000.

As in \cite{Lon2}, denote by
\bea
N_1(\lm, b) &=& \left(\begin{array}{cc}\lm & b\cr
                                0 & \lm \end{array}\right), \qquad {\rm for\;}\lm=\pm 1, \; b\in\R, \lb{3.1}\\
D(\lm) &=& \left(\begin{array}{cc}\lm & 0\cr
                      0 & \lm^{-1} \end{array}\right), \qquad {\rm for\;}\lm\in\R\bs\{0, \pm 1\}, \lb{3.2}\\
R(\th) &=& \left(\begin{array}{cc}\cos\th & -\sin\th \cr
                           \sin\th & \cos\th \end{array}\right), \qquad {\rm for\;}\th\in (0,\pi)\cup (\pi,2\pi), \lb{3.3}\\
N_2(e^{\th\sqrt{-1}}, B) &=& \left(\begin{array}{cc} R(\th) & B \cr
                  0 & R(\th)\end{array}\right), \qquad {\rm for\;}\th\in (0,\pi)\cup (\pi,2\pi)\;\; {\rm and}\; \nn\\
        && \qquad B=\left(\begin{array}{cc} b_1 & b_2\cr
                                  b_3 & b_4\end{array}\right)\; {\rm with}\; b_j\in\R, \;\;
                                         {\rm and}\;\; b_2\not= b_3. \lb{3.4}\eea
Here $N_2(e^{\th\sqrt{-1}}, B)$ is non-trivial if $(b_2-b_3)\sin\theta<0$, and trivial if $(b_2-b_3)\sin\theta>0$.

As in \cite{Lon2}, the $\diamond$-sum (direct sum) of any two block-wise real matrices is defined by
\bea \left(\begin{array}{cc}A_1 & B_1\cr C_1 & D_1\end{array}\right)_{2i\times 2i}\diamond
      \left(\begin{array}{cc}A_2 & B_2\cr C_2 & D_2\end{array}\right)_{2j\times 2j}
=\left(\begin{array}{cccc}A_1 & 0 & B_1 & 0 \cr
                                   0 & A_2 & 0& B_2\cr
                                   C_1 & 0 & D_1 & 0 \cr
                                   0 & C_2 & 0 & D_2\end{array}\right). \nn\eea

For every $M\in\Sp(2n)$, the homotopy set $\Omega(M)$ of $M$ in $\Sp(2n)$ is defined by
$$ \Om(M)=\{N\in\Sp(2n)\,|\,\sg(N)\cap\U=\sg(M)\cap\U\equiv\Gamma\;\mbox{and}
                    \;\nu_{\om}(N)=\nu_{\om}(M)\, \forall\om\in\Gamma\}, $$
where we denote by $\sg(M)$ the spectrum of $M$, by $\Om^0(M)$ the path connected component of $\Om(M)$
containing $M$, and $\nu_{\om}(M)\equiv\dim_{\C}\ker_{\C}(M-\om I)$ for $\om\in\U$.

{\bf Lemma 3.1.} (cf. \cite{Lon2}, Lemma 9.1.5 and List 9.1.12 of \cite{Lon3}) {\it For $M\in\Sp(2n)$ and
$\om\in\U$, the splitting number $S_M^\pm(\om)$ (cf. Definition 9.1.4 of \cite{Lon3}) satisfies
\bea
S_M^{\pm}(\om) &=& 0, \qquad {\it if}\;\;\om\not\in\sg(M).  \lb{3.5}\\
S_{N_1(1,a)}^+(1) &=& \left\{\begin{array}{cc}1, &\quad {\rm if}\;\; a\ge 0, \cr
                                     0, &\quad {\rm if}\;\; a< 0. \end{array}\right. \lb{3.6}\eea

For any $M_i\in\Sp(2n_i)$ with $i=0$ and $1$, there holds }
\be S^{\pm}_{M_0\diamond M_1}(\om) = S^{\pm}_{M_0}(\om) + S^{\pm}_{M_1}(\om),
    \qquad \forall\;\om\in\U. \lb{3.7}\ee

We have the following decomposition theorem

\medskip

{\bf Theorem 3.2.} (cf. \cite{Lon2} and Theorem 1.8.10 of \cite{Lon3}) {\it For any $M\in\Sp(2n)$, there
exists a path $f\in C([0,1],\Om^0(M))$ such that $f(0)=M$ and
\be f(1) = M_1\diamond\cdots\diamond M_k,  \lb{3.8}\ee
where each $M_i$ is a basic normal form listed in (\ref{3.1})-(\ref{3.4}) for $1\leq i\leq k$.}

\medskip
For every $\ga\in\mathcal{P}_\tau(2n)\equiv\{\ga\in C([0,\tau],Sp(2n))\ |\ \ga(0)=I_{2n}\}$, we extend
$\ga(t)$ to $t\in [0,m\tau]$ for every $m\in\N$ by
\be \ga^m(t)=\ga(t-j\tau)\ga(\tau)^j \qquad \forall\;j\tau\le t\le (j+1)\tau \;\;
               {\rm and}\;\;j=0, 1, \ldots, m-1, \lb{3.9}\ee
as in P.114 of \cite{Lon1}. As in \cite{LoZ} and \cite{Lon3}, we denote the Maslov-type indices of
$\ga^m$ by $(i(\ga,m),\nu(\ga,m))$. By \cite{Liu} and \cite{LLo}, the Morse indices of an oriented closed
geodesic $c$ coincide to the Maslov-type indices of the fundamental solution $\ga_c$ of the linearized
Hamiltonian system at $c$, i.e., $(i(c^m),\nu(c^m)) = (i(\ga_c,m),\nu(\ga_c,m))$.

In order to find the iterates which do not change the nullity of a closed geodesic, we distinguish symplectic
matrices with constant or various nullities for all their iterates.
\bea
\Sp_{cnu}(2n) &=& \{M\in \Sp(2n)\;|\; \dim(M^m-I)=\dim(M-I),\;\forall\;m\in\N\},  \nn\\
\Sp_{vnu}(2n) &=& \Sp(2n)\bs \Sp_{cnu}(2n).  \nn\eea
Note that every $M\in \Sp_{vnu}(2n)$ must possess at least one eigenvalue $e^{\sqrt{-1}\th}\in\sg(M)$
with $\th\in (0,2\pi)\cap\pi\Q$. Then for every $M\in\Sp(2n)$ we define
\be  \check{m}(M) = \left\{\begin{array}{cc}
    \min\{k\in\N\;|\;k\th\in 2\pi\N,\;e^{\sqrt{-1}\th}\in\sg(M)\;{\rm with}\;\th\in (0,2\pi)\cap \pi\Q\},
            &{\rm if}\;M\in\Sp_{vnu}(2n), \cr
    +\infty, & {\rm if}\;M\in\Sp_{cnu}(2n). \end{array}\right.  \lb{3.10}\ee
Note that $\check{m}(M)\ge 2$ holds for every $M\in\Sp(2n)$. For a closed geodesic $c$ on a Finsler manifold
$(M,F)$, using the linearized Poincar\'e map $P_c$ we define $\check{m}(c) = \check{m}(P_c)$ correspondingly.
When $(M,F)$ is bumpy, then $\check{m}(c)=+\infty$ holds for every closed geodesic on $(M,F)$.

The following iteration formula and common selection theorem were proved in \cite{LoZ}.

\medskip

{\bf Theorem 3.3.} (cf. Theorem 2.1 and Corollary 2.1 of \cite{LoZ} or Theorem 9.3.1 and Corollary 9.3.2
of \cite{Lon3}) {\it For any path $\ga\in\mathcal{P}_\tau(2n)$, let $M=\ga(\tau)$ and
$C(M)=\sum_{0<\th<2\pi}S_M^-(e^{\sqrt{-1}\th})$. Then for any $m\in\N$ we have
\bea i(\ga,m)
&=& m(i(\ga,1)+S^+_{M}(1)-C(M))\nn\\
& & + 2\sum_{\th\in(0,2\pi)}E\left(\frac{m\th}{2\pi}\right)S^-_{M}(e^{\sqrt{-1}\th}) - (S_M^+(1)+C(M)), \lb{3.11}\eea
and
\be \hat{i}(\ga,1) = i(\ga,1) + S^+_{M}(1) - C(M) + \sum_{\th\in(0,2\pi)}\frac{\th}{\pi}S^-_{M}(e^{\sqrt{-1}\th}).
               \lb{3.12}\ee}

\medskip

{\bf Theorem 3.4.} (cf. Theorem 4.1 of \cite{LoZ} or Theorem 11.1.1 of \cite{Lon3}) {\it Fix an integer $q>0$.
Let $\mu_k\ge 0$ and $\bb_k$ be integers for all $k=1,\cdots,q$. Let $\aa_{k,j}$ be positive numbers for
$j=1,\cdots,\mu_k$ and $k=1,\cdots,q$. Let $\dl\in(0,\frac{1}{2})$ satisfying
$\dl\max_{1\le k\le q}\mu_k<\frac{1}{2}$. Suppose $D_k \equiv \bb_k+\sum_{j=1}^{\mu_k}\aa_{k,j}>0$ for
$k=1,\cdots,q$. Then there exist infinitely many $(N, m_1,\cdots,m_q)\in\N^{q+1}$ such that
\bea
m_k\bb_k+\sum_{j=1}^{\mu_k}E(m_k\aa_{k,j}) = N+\Delta_k,\qquad &&\forall\quad k=1,\cdots,q,   \lb{3.13}\\
\min\{\{m_k\aa_{k,j}\}, 1-\{m_k\aa_{k,j}\}\} < \dl,\qquad &&\forall\ j=1,\cdots,\mu_k,\quad k=1,\cdots,q, \lb{3.14}\\
m_k\aa_{k,j}\in\N,\qquad &&\mbox{if} \qquad \aa_{k,j}\in\Q,   \lb{3.15}\eea
where
\be \Delta_k=\sum_{0<\{m_k\aa_{k,j}\}<\dl}1,\qquad k=1,\cdots,q.\lb{3.16}\ee}

\vskip 1 mm

The common index jump theorem (Theorem 4.3 of \cite{LoZ}) for symplectic paths established by Long
and Zhu in 2002 has become one of the main tools to study the multiplicity and stability problems of
closed solution orbits in Hamiltonian and symplectic dynamics. Here following the ideas of \cite{LoZ},
we prove an enhanced version of it for our usage in this paper by suitably modifying the corresponding
proof in \cite{LoZ} to obtain the index properties of $\ga_k^{2m_k}$ and $\ga_k^{2m_k\pm m}$ with
suitable $m$s.

\medskip

{\bf Theorem 3.5.} ({\bf The enhanced common index jump theorem for symplectic paths}) {\it Let
$\gamma_k\in\mathcal{P}_{\tau_k}(2n)$ for $k=1,\cdots,q$ be a finite collection of symplectic paths.
Let $M_k=\ga_k(\tau_k)$. We extend $\ga_k$ to $[0,+\infty)$ by (\ref{3.9}) inductively. Suppose
\be  \hat{i}(\ga_k,1) > 0, \qquad \forall\ k=1,\cdots,q.  \lb{3.17}\ee
We define
\be  \check{m} \equiv \check{m}(\ga_1,\ldots,\ga_k) = \min\{\check{m}(M_k)\;|\;1\le k\le q\}.  \lb{3.18}\ee
Then for every integer $\bar{m}\in \N$, there exist infinitely many $(q+1)$-tuples
$(N, m_1,\cdots,m_q) \in \N^{q+1}$ such that we have
\be  \nu(\ga_k,2m_k-m) = \nu(\ga_k,2m_k+m) = \nu(\ga_k,1), \quad \forall\;1\le k\le q, 1\le m<\check{m},  \lb{3.19}\ee
and the following hold for all $1\le k\le q$ and $1\le m\le \bar{m}$,
\bea
\nu(\ga_k,2m_k-m) &=& \nu(\ga_k,2m_k+m) = \nu(\ga_k, m),   \lb{3.20}\\
i(\ga_k,2m_k+m) &=& 2N+i(\ga_k,m),                         \lb{3.21}\\
i(\ga_k,2m_k-m) &=& 2N-i(\ga_k,m)-2(S^+_{M_k}(1)+Q_k(m)),  \lb{3.22}\\
i(\ga_k, 2m_k)&=& 2N -(S^+_{M_k}(1)+C(M_k)-2\Delta_k),     \lb{3.23}\eea
where as in (4.45) of \cite{LoZ}, we let
\be \Delta_k = \sum_{0<\{m_k\th/\pi\}<\delta}S^-_{M_k}(e^{\sqrt{-1}\th}), \lb{3.24}\ee
and we define
\be Q_k(m) = \sum_{e^{\sqrt{-1}\th}\in\sg(M_k),\atop \{\frac{m_k\th}{\pi}\}
                   = \{\frac{m\th}{2\pi}\}=0}S^-_{M_k}(e^{\sqrt{-1}\th}). \lb{3.25}\ee}

\medskip

{\bf Proof.} Note first that $\check{m}(q_1,\ldots,q_k)\ge 2$ holds always by the definitions (\ref{3.10})
and (\ref{3.18}). We follow the proof of Theorem 4.2 of \cite{LoZ} and make modifications corresponding to
$1\le m\le \bar{m}$.

Let
\be \dl_0 = \min_{1\le k\le q}\left\{\left.\frac{1}{2},\left\{\frac{h\th}{2\pi}\right\},1-\left\{\frac{h\th}{2\pi}\right\}
   \ \right|\ \frac{\th}{\pi}\in (0,2)\bs\Q, e^{\sqrt{-1}\th}\in\sg(M_k), 1\le h\le \bar{m}\right\}, \qquad \lb{3.26}\ee
and for $1\le k\le q$ and $m\in\N$, let
\bea
&& C(M_k) = \sum_{\th\in(0,2\pi)}S^-_{M_k}(e^{\sqrt{-1}\th}),  \lb{3.27}\\
&& \rho_k = i(\ga_k,1) + S_{M_k}^+(1) - C(M_k),   \lb{3.28}\\
&& I(k,m) = m\rho_k + \sum_{\th\in(0,2\pi)}E\left(\frac{m\th}{\pi}\right)S^-_{M_k}(e^{\sqrt{-1}\th}).  \lb{3.29}\eea
By definition, we have $\dl_0\in (0,1/2]$.

Now we apply Theorem 3.4 to our case. Let $\delta\in(0,\delta_0)$, $\beta_k=\rho_k$, $D_k=\hat{i}(\ga_k,1)$,
$\mu_k=\sum_{\th\in(0,2\pi)}S^-_{M_k}(e^{\sqrt{-1}\th})$, and
$\aa_{k,j} = \frac{\theta_j}{\pi}$, where $e^{\sqrt{-1}\th_j}\in\sigma(M_k)$ for all $1\le j\le \mu_k$ and $1\le k\le q$.
Then for any $\dl\in (0,\dl_0)$, by Theorem 3.4, there exist infinitely many $(q+1)$-tuples $(N,m_1,\cdots,m_q)\in\N^{q+1}$
such that for $1\le k\le q$ we have
\bea
&& \frac{m_k\theta}{\pi} \in \Z, \quad \mbox{whenever}\;\;\frac{\theta}{\pi}\in\Q\cap(0,2)\;\;\mbox{and}\;\;
                     e^{\sqrt{-1}\theta}\in\sigma(M_k),     \lb{3.30}\\
&& \min\left\{\left\{\frac{m_k\th}{\pi}\right\}, 1-\left\{\frac{m_k\th}{\pi}\right\}\right\} <  \dl, \quad  \mbox{whenever}\quad
                     e^{\sqrt{-1}\th}\in\sigma(M_k),     \lb{3.31}\\
&& I(k,m_k) = N+\Delta_k,  \lb{3.32}\eea
where $\Delta_k$ is defined in (\ref{3.24}).

Note that by (\ref{3.30}), the right hand side in (\ref{3.24}) and the left hand side in (\ref{3.31})
are only evaluated on those $e^{\sqrt{-1}\theta}\in\sg(M_k)$ with $\th\in (0,2\pi)\bs \pi\Q$.

\medskip

Here the proof of (\ref{3.19}) is omitted, because it follows directly from that in the Step 2 in the proof
of Theorem 4.3 in \cite{LoZ} by changing $2m_k\th \pm\th$ to $2m_k\th \pm m\th$ for $1\le m<\check{m}$ and
using the definition of $\check{m}$. Next we prove (\ref{3.20})-(\ref{3.23}) in four steps.

\medskip

{\bf Step 1.} {\it Verification of (\ref{3.20})}.

Note first that (\ref{3.20}) holds or not is determined by those eigenvalues $e^{\sqrt{-1}\th}\in\sg(M_k)$
with $\frac{\th}{\pi}\in \Q\cap (0,2)$. For every such eigenvalue $e^{\sqrt{-1}\th}$, by (\ref{3.30})
we have always $2m_k\th \in 2\pi\Z$. Thus $m\th\in 2\pi\Z$ holds if and only if $2m_k\th + m\th\in 2\pi\Z$
holds, if and only if $2m_k\th - m\th\in 2\pi\Z$ holds. Thus at $(2m_k\pm m)$-th iterates, the multiplicity
of the eigenvalue $1$ in $\sg(\ga_k(\tau_k)^{2m_k\pm m})$ is the same as that in $\sg(\ga_k(\tau_k)^{m})$.
This proves (\ref{3.20}).

\medskip

{\bf Step 2.} {\it Verification of (\ref{3.21})}.

\medskip

Now we can compute the value of $i(\ga_k, {2m_k+m})$ with $1\le m\le \bar{m}$ as follows
\bea i(\ga_k,{2m_k+m})
&=& 2m_k(i(\ga_k,1)+S^+_{M_k}(1)-C(M_k)) + 2\sum_{\th\in(0,2\pi)}E\left(\frac{m_k\th}{\pi}\right)S^-_{M_k}(e^{\sqrt{-1}\th}) \nn\\
&&+ 2\sum_{\th\in(0,2\pi)}E\left(\frac{m\th}{2\pi}\right)S^-_{M_k}(e^{\sqrt{-1}\th})-(S^+_{M_k}(1)+C(M_k))\nn\\
&&+ m(i(\ga_k,1)+S^+_{M_k}(1)-C(M_k))+ 2\sum_{\th\in(0,2\pi)}\xi^+_m(m_k,\th)S^-_{M_k}(e^{\sqrt{-1}\th}) \nn\\
&=& 2I(k,m_k)+i(\ga_k,{m}) + 2\sum_{\th\in(0,2\pi)}\xi^+_m(m_k,\th)S^-_{M_k}(e^{\sqrt{-1}\th}) \nn\\
&=& 2(N+\Delta_k) + i(\ga_k,{m}) + 2\sum_{\th\in(0,2\pi)}\xi^+_m(m_k,\th)S^-_{M_k}(e^{\sqrt{-1}\th}) \nn\\
&=& 2N + i(\ga_k,{m}) + 2\sum_{\th\in(0,2\pi)}\xi^+_m(m_k,\th)S^-_{M_k}(e^{\sqrt{-1}\th})+2\Delta_k,
                    \lb{3.33}\eea
where $\xi^+_m(m_k,\th)$ for $1\le m\le \bar{m}$ is defined by
\bea \xi^+_m(m_k,\th)
&\equiv& E\left(\frac{m_k\th}{\pi}+\frac{m\th}{2\pi}\right) - E\left(\frac{m_k\th}{\pi}\right)
       - E\left(\frac{m\th}{2\pi}\right)   \nn\\
&=& E\left(\left\{\frac{m_k\th}{\pi}\right\}+\left\{\frac{m\th}{2\pi}\right\}\right)
       - E\left(\left\{\frac{m_k\th}{\pi}\right\}\right) - E\left(\left\{\frac{m\th}{2\pi}\right\}\right).\lb{3.34}\eea

\medskip

{\bf Claim 1.} {\it For $e^{\sqrt{-1}\th}\in\sg(M_k)$ and any $1\le m\le \bar{m}$, there holds
\be \xi^+_m(m_k,\th)= \left\{\begin{array}{cc}
       0, & \quad {\rm if}\quad \left\{\frac{m_k\th}{\pi}\right\}=0\quad \mbox{or}
                      \quad 1-\delta<\left\{\frac{m_k\th}{\pi}\right\}<1,\cr
      -1, & \quad {\rm if}\quad 0<\left\{\frac{m_k\th}{\pi}\right\}<\delta. \end{array}\right.\lb{3.35}\ee}

In fact, if $\left\{\frac{m_k\th}{\pi}\right\}=0$ holds, $\xi^+_m(m_k,\th)=0$ follows from the definition
(\ref{3.34}) directly.

Note that whenever $\left\{\frac{m_k\th}{\pi}\right\}>0$, then $\th\in (0,2\pi)\bs \pi\Q$ must hold by
(\ref{3.30}) and thus $\{\frac{m\th}{2\pi}\}>0$ for every $1\le m\le \bar{m}$ in the following two subcases.

If $1-\delta<\left\{\frac{m_k\th}{\pi}\right\}<1$, then there holds $1<\left\{\frac{m_k\th}{\pi}\right\}+\delta<2$.
Together with the definition (\ref{3.26}) of $\dl_0$, it yields
$$  1 < \left\{\frac{m_k\th}{\pi}\right\}+\dl < \left\{\frac{m_k\th}{\pi}\right\}+\dl_0
         \le \left\{\frac{m_k\th}{\pi}\right\}+\left\{\frac{m\th}{2\pi}\right\} < 2  $$
for any $1\le m\le \bar{m}$. So we obtain $\xi^+_m(m_k,\th)=0$.

If $0<\left\{\frac{m_k\th}{\pi}\right\}<\dl$, then we have
$$  0 < \left\{\frac{m_k\th}{\pi}\right\}+\left\{\frac{m\th}{2\pi}\right\} < \dl + \left\{\frac{m\th}{2\pi}\right\}
      < \dl_0 + \left\{\frac{m\th}{2\pi}\right\} \le 1   $$
for all $1\le m\le \bar{m}$. Thus by (\ref{3.34}) we get $\xi^+_m(m_k,\th)=-1$. Claim 1 is proved.

\medskip

Now (\ref{3.21}) follows from (\ref{3.24}), (\ref{3.33}) and Claim 1.

\medskip

{\bf Step 3.} {\it Verification of (\ref{3.22}).}

\medskip

Similarly we can compute the value of $i(\ga_k,{2m_k-m})$ with $1\le m\le \bar{m}$ as follows
\bea i(\ga_k,{2m_k-m})
&=& 2m_k(i(\ga_k,1)+S^+_{M_k}(1)-C(M_k)) + 2\sum_{\th\in(0,2\pi)}E\left(\frac{m_k\th}{\pi}\right)S^-_{M_k}(e^{\sqrt{-1}\th}) \nn\\
&& - m(i(\ga_k,1)+S^+_{M_k}(1)-C(M_k)) \nn\\
&& -2\sum_{\th\in(0,2\pi)}E\left(\frac{m\th}{2\pi}\right)S^-_{M_k}(e^{\sqrt{-1}\th})-(S^+_{M_k}(1)+C(M_k)) \nn\\
&& + 2\sum_{\th\in(0,2\pi)}\xi^-_m(m_k,\th)S^-_{M_k}(e^{\sqrt{-1}\th})-2(S^+_{M_k}(1)+C(M_k)) \nn\\
&=& 2I(k,m_k)-i(\ga_k,{m}) + 2\sum_{\th\in(0,2\pi)}\xi^-_m(m_k,\th)S^-_{M_k}(e^{\sqrt{-1}\th})-2(S^+_{M_k}(1)+C(M_k)) \nn\\
&=& 2(N+\Delta_k) - i(\ga_k,{m}) + 2\sum_{\th\in(0,2\pi)}\xi^-_m(m_k,\th)S^-_{M_k}(e^{\sqrt{-1}\th})-2(S^+_{M_k}(1)+C(M_k)) \nn\\
&=& 2N - i(\ga_k,m) - 2S^+_{M_k}(1) \nn\\
&&\quad + 2\sum_{\th\in(0,2\pi)}\xi^-_m(m_k,\th)S^-_{M_k}(e^{\sqrt{-1}\th}) - 2(C(M_k)-\Dl_k),   \lb{3.36}\eea
where $\xi^-_m(m_k,\th)$ for $1\le m\le \bar{m}$ is defined by
\bea \xi^-_m(m_k,\th)
&\equiv& E\left(\frac{m_k\th}{\pi}-\frac{m\th}{2\pi}\right) - E\left(\frac{m_k\th}{\pi}\right)
       + E\left(\frac{m\th}{2\pi}\right)   \nn\\
&=& E\left(\left\{\frac{m_k\th}{\pi}\right\}-\left\{\frac{m\th}{2\pi}\right\}\right)
       - E\left(\left\{\frac{m_k\th}{\pi}\right\}\right) + E\left(\left\{\frac{m\th}{2\pi}\right\}\right).\lb{3.37}\eea

\medskip

{\bf Claim 2.} {\it For $e^{\sqrt{-1}\th}\in\sg(M_k)$ and any $1\le m\le \bar{m}$, there holds }
\be \xi^-_m(m_k,\th)= \left\{\begin{array}{cccc}
              0, & \quad {\rm if}\, \left\{\frac{m_k\th}{\pi}\right\}=0\, {\rm and}\, \left\{\frac{m\th}{2\pi}\right\}=0,\cr
              1, & \quad {\rm if}\, \left\{\frac{m_k\th}{\pi}\right\}=0\, {\rm and}\, \left\{\frac{m\th}{2\pi}\right\}>0,\cr
              1, & \quad {\rm if}\  1-\delta<\left\{\frac{m_k\th}{\pi}\right\}<1, \cr
              0, & \quad {\rm if}\  0<\left\{\frac{m_k\th}{\pi}\right\}<\delta.  \end{array}\right.\lb{3.38}\ee

In fact, if $\left\{\frac{m_k\th}{\pi}\right\}=0$, the two possible values of $\xi^-_m(m_k,\th)$ in terms of that
of $\{\frac{m\th}{2\pi}\}$ follow directly from the definition (\ref{3.37}).

Note also that whenever $\left\{\frac{m_k\th}{\pi}\right\}>0$, then $\th\in (0,2\pi)\bs \pi\Q$ must hold by
(\ref{3.30}) and thus $\{\frac{m\th}{2\pi}\}>0$ for every $1\le m\le \bar{m}$ in the following two subcases.

If $1-\delta<\left\{\frac{m_k\th}{\pi}\right\}<1$, then together with the definition (\ref{3.26}) of $\dl_0$
and the fact $\dl\in (0,\dl_0)$, it yields
$$  1 > \left\{\frac{m_k\th}{\pi}\right\} - \left\{\frac{m\th}{2\pi}\right\} > 1 - 2\dl >0 $$
for all $1\le m\le \bar{m}$. Thus by (\ref{3.37}) we get $\xi^-_m(m_k,\th)=1$.

If $0<\left\{\frac{m_k\th}{\pi}\right\}<\delta$, then $-1<\left\{\frac{m_k\th}{\pi}\right\}-\delta<0$ holds.
Together with the fact $\delta\in(0,\delta_0)$ and (\ref{3.26}), it yields
$$  -1 < \left\{\frac{m_k\th}{\pi}\right\}-\left\{\frac{m\th}{2\pi}\right\} < \left\{\frac{m_k\th}{\pi}\right\}-\dl_0 < 0, \quad
    \forall\; 1\le m\le \bar{m}.  $$
So by (\ref{3.37}) we obtain $\xi^-_m(m_k,\th)=0$. Then Claim 2 is proved.

Now because for each $k$ and $m$, the integers $\Dl_k$ and $Q_k(m)$ defined in (\ref{3.24}) and (\ref{3.25})
just count the multiplicities of those eigenvalues $e^{\sqrt{-1}\th}\in\sg(M_k)$ which contribute nothing
to $\xi^-_m(m_k,\th)$ in the forth and the first cases of (\ref{3.38}) respectively, Claim 2 yields
\bea \sum_{\th\in(0,2\pi)}\xi^-_m(m_k,\th)S^-_{M_k}(e^{\sqrt{-1}\th}) = C(M_k) - \Dl_k- Q_k(m)  \lb{3.39}\eea
for $1\le k\le q$ and $1\le m\le \bar{m}$. Together with (\ref{3.36}) it yields
\bea i(\ga_k,{2m_k-m})
&=& 2N - i(\ga_k,m) - 2(S^+_{M_k}(1)+C(M_k)-\Dl_k) + 2(C(M_k)-\Dl_k-Q_k(m))  \nn\\
&=& 2N - i(\ga_k,m) - 2(S^+_{M_k}(1)+Q_k(m)), \lb{3.40}\eea
i.e., (\ref{3.22}) holds.

\medskip

{\bf Step 4.} {\it Verification of (\ref{3.23}).}

\medskip

By Theorem 3.3, (\ref{3.29}) and (\ref{3.32}) we have
\bea i(\ga_k,{2m_k})
&=& 2m_k(i(\ga_k,1)+S^+_{M_k}(1)-C(M_k))
    + 2\sum_{\th\in(0,2\pi)}E\left(\frac{m_k\th}{\pi}\right)S^-_{M_k}(e^{\sqrt{-1}\th}) \nn\\
& & \qquad - (S^+_{M_k}(1)+C(M_k))          \nn\\
&=& 2I(k,m_k) - (S^+_{M_k}(1)+C(M_k)) \nn\\
&=& 2(N+\Dl_k) - (S^+_{M_k}(1)+C(M_k)) \nn\\
&=& 2N - (S^+_{M_k}(1)+C(M_k)-2\Dl_k), \lb{3.41}\eea
i.e., (\ref{3.23}) holds.

The proof of Theorem 3.5 is complete. \hfill\hb

\medskip

{\bf Remark 3.6.} (i) Setting $\bar{m}=1$, our Theorem 3.5 yields the common index jump Theorem 4.3 of \cite{LoZ}
(cf. Theorem 11.2.1 in \cite{Lon3}). As we mentioned in Section 1, Theorem 3.5 allows to get precise index information
of every iterates $\ga_k^{2m_k\pm m}$ for $1\le m\le \bar{m}$ with any given $\bar{m}>0$. This is the first improvement
of CIJT. The second improvement is that we get the precise index information of $\ga_k^{2m_k}$ in (\ref{3.23}) instead
of given only the lower and upper bounds on it as in Theorem 4.3 of \cite{LoZ}. These two improvements play important
and crucial roles in proofs of our main theorems on closed geodesics in Section 4 below.

(ii) By (4.10) and (4.40) in \cite{LoZ} (cf. (11.1.10) and (11.2.14) of \cite{Lon3}), we have
\be  m_k=\left(\left[\frac{N}{\bar{M}\hat i(\gamma_k,1)}\right]+\chi_k\right)\bar{M},\quad 1\le k\le q,\lb{3.42}\ee
where $\chi_k=0$ or $1$ for $1\le k\le q$, $\bar{M}$ is a positive integer such that $\frac{\bar{M}\theta}{\pi}\in\Z$
whenever $e^{\sqrt{-1}\theta}\in\sigma(M_k)$ and $\frac{\theta}{\pi}\in\Q$ for some $1\le k\le q$, and we set
$\bar{M}=1$ if no such eigenvalues exist. Furthermore, by (4.20) in Theorem 4.1 of \cite{LoZ} (cf. (11.1.20) of
\cite{Lon3}), for any $\ep>0$, we can choose $N$ and $\{\chi_k\}_{1\le k\le q}$ such that
\be  \left|\left\{\frac{N}{\bar{M}\hat{i}(\ga_k,1)}\right\}-\chi_k\right|<\ep,\qquad 1\le k\le q.  \lb{3.43}\ee

(iii) Let  $\mu_k=\sum_{\theta\in(0,2\pi)}S_{M_k}^-(e^{\sqrt{-1}\theta})$ for $1\le k\le q$, and
$\alpha_{k, j}=\frac{\theta_j}{\pi}$ where $e^{\sqrt{-1}\theta_j}\in\sigma(M_k)$ for $1\le j\le \mu_k$ and
$1\le k\le q$.  Let $l=q+\sum_{1\le k\le q}\mu_k$ and
\be
v = \left(\frac{1}{\bar{M}\hat i(\ga_1, 1)},\ldots, \frac{1}{\bar{M}\hat i(\ga_q, 1)}, \frac{\aa_{1, 1}}{\hat i(\ga_1, 1)},
\frac{\aa_{1, 2}}{\hat i(\ga_1, 1)}, \ldots\frac{\aa_{1,\mu_1}}{\hat i(\ga_1, 1)},
\frac{\aa_{2, 1}}{\hat i(\ga_2, 1)},\ldots, \frac{\aa_{q, \mu_q}}{\hat i(\ga_q,1)}\right)\in\R^l. \lb{3.44}\ee
Then the existence of the $(q+1)$-tuple $(N, m_1,\ldots,m_q) \in \N^{q+1}$ in the above Theorem 3.5 is
equivalent to the existence of a vertex
\be \chi=(\chi_1,\ldots,\chi_q,\chi_{1, 1},\chi_{1, 2},\ldots,\chi_{1, \mu_1},
       \chi_{2, 1},\ldots,\chi_{q, \mu_q})\lb{3.45}\ee
of the unit cube $[0,\,1]^l$ and infinitely many $N\in\N$ such that for any prescribed sufficiently small $\dl>0$,
\be  |\{Nv\}-\chi|<\delta  \lb{3.46}\ee
holds (cf. Pages 346 and 349 of \cite{LoZ}), where $\{mv\}=(\{mv_1\},\ldots,\{mv_l\})$ for $v=(v_1,\ldots,v_l)$
and $m\in\N$. Note that here $\chi$ is some vertex of $[0,1]^l$ belonging to $H\cap[0,1]^l$, where $H$ is the
closure of the set $\{ \{mv\} \,|\, m\in\N\}$ in $[0,1]^l$.

(iv) Note that, given $M_0\in\N$, we can require $N$ to be a multiple of $M_0$, since the closure of the set
$\{\{Nv\}: N\in\N, \;M_0|N\}$ is still a closed additive subgroup of the torus ${\bf T}^h$ for some $h\in\N$,
where $v$ is defined in the above (\ref{3.44}). Then we can use the proof of Step 2 in Theorem 4.1 of \cite{LoZ}
to get $N$.

(v) Very recently Gutt and Kang obtained also independently a variant (Theorem 2.2 of \cite{GuK}) of the common
index jump theorem of Long and Zhu, where they considered only non-degenerate symplectic paths, gave lower and
upper bound index estimates on $\ga_k^{2m_k}$, and studied the indices of iterates $\ga_k^{2m_k\pm m}$ for
$1\le m\le \tilde{m}$ with $\tilde{m}$ being given. Their result is in fact weaker than our above Theorem 3.5 and yields rougher
information. Based on this result, Gutt and Kang obtained an estimate on the number of closed characteristics on
compact star-shaped hypersurfaces in $\R^{2n}$ by assuming that all the closed characteristics and their iterates
are non-degenerate and their minimal Conley-Zehnder indices are at least $n-1$.

\medskip

\setcounter{figure}{0}
\setcounter{equation}{0}
\section{Studies on closed geodesics}

Now applying our enhanced common index jump theorem, we give proofs of the main theorems described in Section 1.

\subsection{Proof of Theorem 1.1}

In order to prove Theorem 1.1, let $(M,F)$ be a compact simply-connected manifold with a bumpy,
irreversible Finsler metric $F$ and satisfy $H^*(M;\Q)\cong T_{d,n+1}(x)$ for some even integer
$d\ge 2$ and some integer $n\ge 1$. We make the following assumption,

\medskip

{\bf (FCG)} {\it Suppose that there exist only finitely many prime closed geodesics
$\{c_k\}_{k=1}^q$ with $i(c_k)\ge 1$ for $k=1,\cdots,q$ on $(M,F)$.}

\medskip

When the Finsler metric $F$ is bumpy, for every prime closed geodesic $c$, by Theorem 3.2 the
basic normal form decomposition of the linearized Poincar\'{e} map $P_{c}$  possesses the following
form
\bea f_c(1)
&=& R(\th_1)\,\dm\,\cdots\,\dm\,R(\th_r)\,\dm\,D(\lm_1)\dm \,\cdots\,\dm D(\lm_s) \nn\\
& & \dm\,N_2(e^{\aa_{1}\sqrt{-1}},A_{1})\,\dm\,\cdots\,\dm\,N_2(e^{\aa_{r_{\ast}}\sqrt{-1}},A_{r_{\ast}})\,
    \dm\,N_2(e^{\bb_{1}\sqrt{-1}},B_{1})\,\dm\,\cdots\,\dm\,N_2(e^{\bb_{r_{0}}\sqrt{-1}},B_{r_{0}}), \nn\eea
where $\lm_j\in \R\bs\{0,\pm 1\}$ for $1\le j\le s$, $\frac{\th_{j}}{2\pi}\not\in\Q$ for $1\le j\le r$,
$\frac{\aa_{j}}{2\pi}\not\in\Q$ for $1\le j\le r_{\ast}$, $\frac{\bb_{j}}{2\pi}\not\in\Q$ for $1\le j\le r_0$,
and
\be r + s + 2r_{\ast} + 2r_0 = dn - 1. \lb{4.1}\ee
Then as proved in \cite{DuL1}, we obtain the index iteration formula of $c^m$
\be i(c^m) = m(i(c)-r) + 2\sum_{j=1}^r\left[\frac{m\th_j}{2\pi}\right] +r,\quad \nu(c^m)=0,
    \qquad\forall\ m\ge 1. \lb{4.2}\ee

\medskip

{\bf Proof of Theorem 1.1.}

\medskip

We carry out the proof in three steps.

\medskip

{\bf Step 1}. {\it The existence of $\frac{dn(n+1)}{4}$ distinct prime closed geodesics.}

\medskip

Since by the assumption (FCG), there exist only finitely many distinct prime closed geodesics
on the bumpy manifold $(M,F)$, any closed geodesic $c_k$ among $\{c_k\}_{k=1}^q$ must have
positive mean index (cf. Theorem 3 of \cite{BaK} or Lemma 3.2 of \cite{DuL2}), i.e.,
\be \hat{i}(c_k)>0,\qquad 1\le k\le q, \lb{4.3}\ee
which implies that $i(c_k^m)\rightarrow +\infty$ as $m\rightarrow +\infty$. So the positive
integer $\bar{m}$ defined by
\be \bar{m}=\max_{1\le k\le q}\left\{\min\{m\in\N\ |\ i(c_k^{m})\ge i(c_k)+2(dn-1)\}\right\}
              \lb{4.4}\ee
is well-defined and finite.

For the integer $\bar{m}$ defined in (\ref{4.4}), by (\ref{4.3}) it follows from Theorem 3.5
that there exist infinitely many $q+1$-tuples $(N, m_1, \ldots, m_q)\in\N^{q+1}$ such that
for any $1\le k\le q$, there holds
\bea
i(c_k^{2m_k-m}) &=& 2N-i(c_k^m),\quad 1\le m\le\bar{m}, \lb{4.5}\\
i(c_k^{2m_k}) &=& 2N-C(M_k)+2\Delta_k,  \lb{4.6}\\
i(c_k^{2m_k+m}) &=& 2N+i(c_k^m),\quad 1\le m\le\bar{m}, \lb{4.7}\eea
where $M_k=P_{c_k}\in \Sp(2(dn-1))$ is the linearized Poincar\'e map of the prime closed
geodesic $c_k$. Here note that in the bumpy case, $S^+_{M_k}(1)=0$ and $Q_k(m)=0$ holds for
all $m\in\N$.

On one hand, there holds $i(c_k^m)\ge i(c_k)$ for any $m\ge 1$ by the Bott-type formulae (cf.
\cite{Bot} and Theorem 9.2.1 of \cite{Lon3}). Thus by (\ref{4.5}), (\ref{4.6}) and the
assumption $i(c_k)\ge 1$ for $1\le k\le q$ in the theorem, it yields
\bea
i(c_k^{2m_k-m}) &=& 2N-i(c_k^m)\le 2N-i(c_k)\le 2N-1,\quad 1\le m\le\bar{m}, \lb{4.8}\\
i(c_k^{2m_k+m}) &=& 2N+i(c_k^m)\ge 2N+i(c_k)\ge 2N+1,\quad 1\le m\le\bar{m}. \lb{4.9}\eea

On the other hand, note that the equalities (\ref{3.33}) and (\ref{3.36}) in the proofs of
Steps 2 and 3 of Theorem 3.5 in fact hold for every $m\in\N$ without using any information
on $\bar{m}$. So by (\ref{3.33}) and the definition (\ref{4.4}) of $\bar{m}$ we obtain
\bea i(c_k^{2m_k+m})
&=& 2N + i(c_k^{m}) + 2\sum_{\th\in(0,2\pi)}\xi^+_m(m_k,\th)S^-_{M_k}(e^{\sqrt{-1}\th})+2\Delta_k \nn\\
&\ge& 2N+i(c_k)+2(dn-1)+ 2\sum_{\th\in(0,2\pi)}\xi^+_m(m_k,\th)S^-_{M_k}(e^{\sqrt{-1}\th})+2\Delta_k \nn\\
&\ge& 2N+i(c_k), \qquad \forall\ m\ge \bar{m},\lb{4.10}\eea
where the first inequality follows from the definition (\ref{4.4}) of $\bar{m}$, and the last
inequality follows from the facts $\xi^+_m(m_k,\th)\in\{-1,0\}$ by Claim 1 and $\Delta_k\ge 0$.

Similarly by (\ref{3.36}) and the definition (\ref{4.4}) of $\bar{m}$ we obtain
\bea i(c_k^{2m_k-m})
&=& 2N - i(c_k^{m})+ 2\sum_{\th\in(0,2\pi)}\xi^-_m(m_k,\th)S^-_{M_k}(e^{\sqrt{-1}\th})-2(S^+_{M_k}(1)+C(M_k)-\Delta_k) \nn\\
&\le& 2N-i(c_k)-2(dn-1)+ 2\sum_{\th\in(0,2\pi)}\xi^-_m(m_k,\th)S^-_{M_k}(e^{\sqrt{-1}\th})-2(C(M_k)-\Delta_k) \nn\\
&\le & 2N-i(c_k), \qquad \forall\ \bar{m}\le m<2m_k, \lb{4.11}\eea
where the first inequality follows from the definition (\ref{4.4}) of $\bar{m}$, and the last
inequality follows from the facts the total algebraic multiplicity of eigenvalues of $M_k$ on $\U$
is at most $2(\dim M -1) = 2(dn-1)$ and $C(M_k)\ge\Delta_k$.

In summary, by (\ref{4.5})-(\ref{4.11}), for $1\le k\le q$, we have proved
\bea
&& i(c_k^{2m_k-m})\le 2N-1,\qquad \forall\ 1\le m\le 2m_k, \lb{4.12}\\
&& i(c_k^{2m_k})=2N-C(M_k)+2\Delta_k,\lb{4.13}\\
&& i(c_k^{2m_k+m})\ge 2N+1,\qquad \forall\ m\ge 1.\lb{4.14}\eea

{\bf Claim 3.} {\it For $N\in\N$ in Theorem 3.5 satisfying (\ref{4.12})-(\ref{4.14}) and $2NB(d,n)\in2\N$,
we have
\be \sum_{1\le k\le q} 2m_k\gamma_{c_k}=2NB(d,n). \lb{4.15}\ee}

In fact, here we follow some ideas in \cite{Wan1} and \cite{Wan2} to prove this Claim. By Theorem 2.3,
(\ref{3.42}) and (\ref{3.43}) of Remark 3.6 with
$$ \ep < \frac{1}{1+2\bar{M}\sum_{1\le k\le q}|\ga_{c_k}|},  $$
it yields
\bea \left|2NB(d,n)-\sum_{k=1}^q 2m_k\ga_{c_k}\right|
&=& \left|\sum_{k=1}^q\frac{2N\ga_{c_k}}{\hat{i}(c_k)}-\sum_{k=1}^q 2\ga_{c_k}
                    \left(\left[\frac{N}{\bar{M}\hat{i}(c_k)}\right]+\chi_k\right)\bar{M}\right| \nn\\
&\le& 2\bar{M}\sum_{k=1}^q |\ga_{c_k}|\left|\left\{\frac{N}{\bar{M}\hat{i}(c_k)}\right\}-\chi_k\right| \nn\\
&<& 2\bar{M}\ep\sum_{k=1}^q|\ga_{c_k}| \nn\\
&<& 1. \lb{4.16}\eea
Since each $2m_k\ga_{c_k}$ is an integer, Claim 3 is proved.

\medskip

Now by Lemma 2.1, Definition 2.2 and Theorem 3.3, it yields
\bea \sum_{m=1}^{2m_k} (-1)^{i(c_k^m)} \dim \ol{C}_{i(c_k^{m})}(E,c_k^m)
&=& \sum_{i=0}^{m_k-1} \sum_{m=2i+1}^{2i+2} (-1)^{i(c_k^m)} \dim \ol{C}_{i(c_k^{m})}(E,c_k^m) \nn\\
&=& \sum_{i=0}^{m_k-1} \sum_{m=1}^{2} (-1)^{i(c_k^m)} \dim \ol{C}_{i(c_k^{m})}(E,c_k^m) \nn\\
&=& m_k \sum_{m=1}^{2} (-1)^{i(c_k^m)} \dim \ol{C}_{i(c_k^{m})}(E,c_k^m) \nn\\
&=& 2m_k\ga_{c_k},\qquad \forall\ 1\le k\le q, \lb{4.17}\eea
where the second equality follows from Lemma 2.1 and the fact $i(c_k^{m+2})-i(c_k^m)\in 2\Z$ for all
$m\ge 1$ from Theorem 3.3, and the last equality follows from Lemma 2.1 and Definition 2.2.

By (\ref{4.14}) and Lemma 2.1, we know that all $c_k^{2m_k+m}$'s with $m\ge1$ and $1\le k\le q$
have no contributions to the alternative sum $\sum_{p=0}^{2N}(-1)^p M_p$. Similarly again by Lemma 2.1
and (\ref{4.12}), all $c_k^{2m_k-m}$'s with $1\le m<2m_k$ and $1\le k\le q$ only have contributions
to $\sum_{p=0}^{2N}(-1)^p M_p$.

Thus for the Morse-type numbers $M_p$'s defined by (\ref{2.10}), by (\ref{4.17}) we have
\bea \sum_{p=0}^{2N}(-1)^p M_p
&=& \sum_{k=1}^{q}\ \sum_{1\le m\le 2m_k \atop i(c_k^{m})\le 2N} (-1)^{i(c_k^m)} \dim \ol{C}_{i(c_k^{m})}(E,c_k^m) \nn\\
&=& \sum_{k=1}^{q}\ \sum_{m=1}^{2m_k} (-1)^{i(c_k^m)}\dim\ol{C}_{i(c_k^{m})}(E,c_k^m) \nn\\
& & -\sum_{1\le k\le q \atop i(c_k^{2m_k})\ge 2N+1} (-1)^{i(c_k^{2m_k})} \dim \ol{C}_{i(c_k^{2m_k})}(E,c_k^{2m_k}) \nn\\
&=& \sum_{k=1}^{q} 2m_k\ga_{c_k}
  -\sum_{1\le k\le q \atop i(c_k^{2m_k})\ge 2N+1} (-1)^{i(c_k^{2m_k})} \dim \ol{C}_{i(c_k^{2m_k})}(E,c_k^{2m_k}).
      \lb{4.18}\eea

In order to exactly know whether the iterate $c_k^{2m_k}$ of $c_k$ has contributions to the alternative
sum $\sum_{p=0}^{2N}(-1)^p M_p(k)$ for $1\le k\le q$, we let
\bea
N_+^e &=& ^{\#}\{1\le k\le q\ |\ i(c_k^{2m_k})\ge 2N+1,\ i(c_k^{2m_k})-i(c_k)\in 2\N_0,\ i(c_k)\in 2\N\}, \lb{4.19}\\
N_+^o &=& ^{\#}\{1\le k\le q\ |\ i(c_k^{2m_k})\ge 2N+1,\ i(c_k^{2m_k})-i(c_k)\in 2\N_0,\ i(c_k)\in 2\N-1\}, \lb{4.20}\\
N_-^e &=& ^{\#}\{1\le k\le q\ |\ i(c_k^{2m_k})\le 2N-1,\ i(c_k^{2m_k})-i(c_k)\in 2\N_0,\ i(c_k)\in 2\N\}, \lb{4.21}\\
N_-^o &=& ^{\#}\{1\le k\le q\ |\ i(c_k^{2m_k})\le 2N-1,\ i(c_k^{2m_k})-i(c_k)\in 2\N_0,\ i(c_k)\in 2\N-1\}. \lb{4.22}\eea
Here note that by (ii) of Remark 3.6, we can suppose that $N$ is a multiple of $D=d(n+1)-2$.

Thus by Theorem 2.3, Claim 3, (\ref{4.18}), the definitions of $N^{e}_+$ and $N_+^o$ and Lemma 2.4, we have
\bea -\frac{Ndn(n+1)}{D}+N_+^o-N_+^e
&=& 2NB(d,n)+N_+^o-N_+^e  \nn\\
&=& \sum_{k=1}^q 2m_k\gamma_{c_k}+N_+^o-N_+^e \nn\\
&=& \sum_{p=0}^{2N}(-1)^p M_p \nn\\
&\ge&\sum_{p=0}^{2N}(-1)^p b_p(\Lm M) \nn\\
&=& -\frac{dn(n+1)}{2D}(2N-d)+\frac{dn(n-1)}{4}-1-\ep_{d,n}(2N-1) \nn\\
&=& -\frac{Ndn(n+1)}{D}+\tilde{N}, \lb{4.23}\eea
where recall $D=d(n+1)-2$ defined in Lemma 2.4. Then it implies
\be N_+^o\ge \tilde{N}\equiv\frac{d^2n(n+1)}{2D}+\frac{dn(n-1)}{4}-1-\ep_{d,n}(2N-1).\lb{4.24}\ee

Note that $N$ is a multiple of $D$. So there holds $\{\frac{2N-1-(d-1)}{D}\}=1-\frac{d}{D}=\frac{dn-2}{D}$.
Then by (\ref{2.9}) we have
\bea \ep_{d,n}(2N-1)
&=& \frac{dn-2}{dn}-\frac{2n+d-2}{dn}\left(1-\frac{d}{D}\right)
               -n\left\{\frac{dn-2}{2}\right\}-\left\{\frac{dn-2}{d}\right\} \nn\\
&=& \frac{(d-2)D+2d-d^2}{Dd}-\left\{-\frac{2}{d}\right\} \nn\\
&=& 1-\frac{2}{d}-\left\{-\frac{2}{d}\right\}-\frac{d-2}{D} \;=\; -\frac{d-2}{D}, \nn\eea
which, together with (\ref{4.24}), yields
\be N_+^o\ge\tilde{N} = \frac{d^2n(n+1)}{2D}+\frac{dn(n-1)}{4}-1+\frac{d-2}{D} = \frac{dn(n+1)}{4}. \lb{4.25}\ee

\medskip

{\bf Step 2.} {\it The existence of another set of $\frac{dn(n+1)}{4}$ distinct prime closed geodesics.}

\medskip

Now in order to understand $\Delta_k$s further, we need to understand a symmetry found in the proof of
the common index jump theorem, i.e., the item (c) of Theorem 4.2 of \cite{LoZ} (cf. Theorem 11.1.2 on
pp.234-235 of \cite{Lon3}, where $\R^n$ there becomes $\R^l$ in Remark 3.6 now). There $A(v)$ is a
slanted sub-space of $\R^l$, which possesses the symmetry $A(v)=-A(v)$. Let $\pi$ be the projection from
$\R^l$ to the torus $\T^l$. Then the pull back under $\pi^{-1}$ of the closure $\hat{A}(v)$ of the
set $\{\{mv\}\;|\;m\in\N\}$ can be viewed as a union of may be more than one slanted sub-cubs of the full
dimensional cub $[0,1]^l$. As mentioned in Remark 3.6, the $(q+1)$-tuple $(N, m_1, \ldots, m_q)$ in CIJT
and Theorem 3.5 is found corresponding to a common vertex of $\hat{A}(v)$ and $[0,1]^l$. By the symmetry
$A(v)=-A(v)$, the vertex $\hat{\chi}=\hat{1}-\chi$ of $[0,1]^l$ opposite to $\chi$ is also in $\hat{A}(v)$
as proved in Lemma 3.18 of \cite{DuL3}, where we denote by $\hat{1}=(1, \ldots, 1)\in\R^l$. This is
illustrated in Figure 4.1 in which the simplest case of $\hat{A}(v)$ with only one connected component
is exhibited. A more complicated case of $\hat{A}(v)$ with more than one connected components is
illustrated in Figure 3.1 of \cite{DuL3}.

\begin{figure}
\begin{center}
\includegraphics[width=8cm,height=8cm]{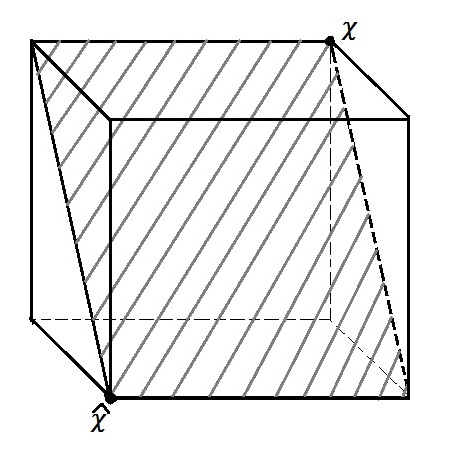}
\vskip -5 mm
\caption{The vertex $\chi$ and the opposite vertex $\hat{\chi}=\hat{1}-\chi$ in $\hat{A}(v)$.}
\end{center}
\end{figure}

Now by Corollary 3.19 of \cite{DuL3} and the above (\ref{4.3}), we can apply Theorem 3.5 again and find
a new $(q+1)$-tuples $(N', m_1', \ldots, m_q')\in\N^{q+1}$ corresponding to the vertex $\hat{\chi}$ such
that similarly to (\ref{4.12})-(\ref{4.14}) for every $1\le k\le q$, there holds
\bea
&& i(c_k^{2m_k'-m})\le 2N'-1,\qquad \forall\ 1\le m<2m_k, \lb{4.26}\\
&& i(c_k^{2m_k'})=2N'-C(M_k)+2\Delta'_k,  \lb{4.27}\\
&& i(c_k^{2m_k'+m})\ge 2N'+1,\qquad \forall\ m\ge 1,  \lb{4.28}\eea
where according to the choice of $\hat{\chi}=\hat{1}-\chi$ symmetric to $\chi$, we have
$$  \Delta'_k = \sum_{0<\{m'_k\th/\pi\}<\delta}S^-_{M_k}(e^{\sqrt{-1}\th}).  $$
Then we have the following

{\bf Claim 4.} {\it For every $1\le k\le q$, }
\be  \Delta_k + \Delta'_k = C(M_k). \lb{4.29}\ee

\medskip

In fact, since the Finsler metric $F$ is bumpy, every eigenvalue $e^{\sqrt{-1}\th_{k_j}}\in\sg(M_k)$
appearing in (\ref{3.11}) must satisfy $\frac{\th_{k_j}}{2\pi}\not\in\Q$ for all $1\le j\le C(M_k)$ and
$1\le k\le q$.

For $\dl\in (0,\dl_0)$ with $\dl_0>0$ small enough chosen in (\ref{3.26}), it follows from (\ref{3.31})
that $\{\frac{m_k\th_{k_j}}{\pi}\}$ must satisfy either $0<\{\frac{m_k\th_{k_j}}{\pi}\}<\delta$ or
$1-\dl < \{\frac{m_k\th_{k_j}}{\pi}\}< 1$.

Since the $(q+1)$-tuple $(N', m_1', \ldots, m_q')\in\N^{q+1}$ is chosen according to the
point $\hat{\chi}$ with some small enough constant $\dl'\in (0,\dl_o)$. Therefore by the symmetry of $\chi$
and $\hat{\chi}=\hat{1}-\chi$ and Corollary 3.19 of \cite{DuL3}, we obtain
\bea
1-\dl' < \{\frac{m_k'\th_{k_j}}{\pi}\}<1 & {\rm if\;and\;only\;if} & 0<\{\frac{m_k\th_{k_j}}{\pi}\}<\dl, \nn\\
0<\{\frac{m'_k\th_{k_j}}{\pi}\}<\dl' & {\rm if\;and\;only\;if} & 1-\dl<\{\frac{m_k\th_{k_j}}{\pi}\}<1. \nn\eea
Then they yield
$$  \Delta'_k = \sum_{0<\{m'_k\th/\pi\}<\dl'}S^-_{M_k}(e^{\sqrt{-1}\th})
       = \sum_{0<1-\{m_k\th/\pi\}<\dl}S^-_{M_k}(e^{\sqrt{-1}\th}),  $$
where the property (\ref{3.31}) is used again whenever $S^-_{M_k}(e^{\sqrt{-1}\th}) >0$. Thus by (\ref{3.31}) and
the definition (\ref{3.24}) of $\Delta_k$ we obtain
\bea \Delta_k + \Delta'_k
&=& \sum_{0<\{m_k\th/\pi\}<\dl}S^-_{M_k}(e^{\sqrt{-1}\th}) + \sum_{1-\dl<\{m_k\th/\pi\}<1}S^-_{M_k}(e^{\sqrt{-1}\th}) \nn\\
&=& \sum_{\th\in (0,2\pi)}S^-_{M_k}(e^{\sqrt{-1}\th}) \;=\; C(M_k), \qquad\forall\ 1\le k\le q, \nn\eea
i.e., Claim 4 holds.

\medskip

Similarly, we define
\bea
N_+^{'e} &=& ^{\#}\{1\le k\le q\ |\ i(c_k^{2m_k'})\ge 2N'+1,\ i(c_k^{2m_k'})-i(c_k)\in 2\N_0,\ i(c_k)\in 2\N\},   \lb{4.30}\\
N_+^{'o} &=& ^{\#}\{1\le k\le q\ |\ i(c_k^{2m_k'})\ge 2N'+1,\ i(c_k^{2m_k'})-i(c_k)\in 2\N_0,\ i(c_k)\in 2\N-1\}, \lb{4.31}\\
N_-^{'e} &=& ^{\#}\{1\le k\le q\ |\ i(c_k^{2m_k'})\le 2N'-1,\ i(c_k^{2m_k'})-i(c_k)\in 2\N_0,\ i(c_k)\in 2\N\},   \lb{4.32}\\
N_-^{'o} &=& ^{\#}\{1\le k\le q\ |\ i(c_k^{2m_k'})\le 2N'-1,\ i(c_k^{2m_k'})-i(c_k)\in 2\N_0,\ i(c_k)\in 2\N-1\}. \lb{4.33}\eea

So by (\ref{4.27}) and (\ref{4.29}) it yields
\be i(c_k^{2m_k'}) = 2N'-C(M_k)+2(C(M_k)-\Delta_k)=2N'+C(M_k)-2\Delta_k. \lb{4.34}\ee

So together with definitions (\ref{4.19})-(\ref{4.22}) and (\ref{4.30})-(\ref{4.33}) it yields
\be N_{\pm}^e = N_{\mp}^{'e}, \qquad N_{\pm}^o = N_{\mp}^{'o}. \lb{4.35}\ee

By Lemma 2.1 and (\ref{4.28}), we know that all $c_k^{2m_k'+m}$'s with $m\ge1$ and $1\le k\le q$ have
no contribution to the alternative sum $\sum_{p=0}^{2N'}(-1)^p M_p$. Similarly also by Lemma 2.1 and
(\ref{4.26}), all $c_k^{2m_k'-m}$'s with $m\ge1$ and $1\le k\le q$ only have contribution to
$\sum_{p=0}^{2N'}(-1)^p M_p$.

Thus, through carrying out arguments similar to (\ref{4.23})-(\ref{4.25}), by Claim 3, the
definitions of $N^{'e}_+$ and $N^{'o}_+$ and Lemma 2.4, together with (\ref{4.35}), we obtain
\be N_-^{o} = N_+^{'o}\ge\tilde{N'} = \frac{dn(n+1)}{4}. \lb{4.36}\ee

So by (\ref{4.23}) and (\ref{4.36}) we get
\be q \ge N_+^{o}+N_-^{o} \ge \frac{dn(n+1)}{4}+\frac{dn(n+1)}{4} = \frac{dn(n+1)}{2}. \lb{4.37}\ee

{\bf Step 3.} In addition, any hyperbolic closed geodesic $c_k$ must have $i(c_k^{2m_k})=2N$ since
there holds $C(M_k)=0$ in the hyperbolic case. However, by (\ref{4.20}) and (\ref{4.22}), there exist
at least $(N_+^{o}+N_-^{o})$ prime closed geodesics with odd indices $i(c_k^{2m_k})$. So all these
$(N_+^{o}+N_-^{o})$ closed geodesics are non-hyperbolic. Then (\ref{4.38}) shows that there exist
at least $\frac{dn(n+1)}{2}$ distinct non-hyperbolic closed geodesics. And (\ref{4.20}), (\ref{4.22})
and (\ref{4.37}) show that all these non-hyperbolic closed geodesics and their iterations have odd
Morse indices. This completes the proof of Theorem 1.1. \hfill\hb

\medskip

\subsection{Proofs of other main theorems}

{\bf Proof of Theorem 1.3.}

\medskip

In fact, let $(M,F)$ be a compact,simply-connected, bumpy Finsler manifold satisfying $\dim M \in 2\N$. We suppose
$\bar{K}_-(M,F)\ge 0$. Then the first conclusion in Theorem 1.3 follows from the proof of Synge Theorem in Finsler
geometry, cf. Theorem 8.8.1 and its proof in Pages 221-223 of \cite{BCS}. For the readers'
convenience, we sketch its proof here.

Note first that because $M$ is compact, simply-connected, there exists at least one prime closed geodesic
on $(M,F)$. As in our discussion near Definition 1.2, let $c$ be any prime closed geodesic on $(M,F)$
satisfying $F(\dot{c}(t))=1$, then $U_{\dot{c},u}$ is well defined and $\bar{K}_-(M,F)<+\infty$ holds.
Then by the second variation formula of the energy functional $E$ at $c$ (cf. Section 2 of \cite{Rad3}) and
the assumption $\bar{K}_-(M,F)\ge 0$, we have
\be  E^{\prime\prime}(c)(U_{\dot{c},u}, U_{\dot{c},u})
          = -\int_0^1K(\dot{c}, U_{\dot{c},u})dt \le -\bar{K}_-(M,F) \le 0. \lb{4.38}\ee
On the other hand, since $F$ is bumpy, then $c$ is non-degenerate. Thus we have $\nu(c)=0$.
If $i(c)=0$, then $E^{\prime\prime}(c)$ must be strictly positive definite, and specially we have
$E^{\prime\prime}(c)(U_{\dot{c},u}, U_{\dot{c},u})>0$. It contradicts to (\ref{4.38}) and proves
the first conclusion in Theorem 1.3.

Then the second conclusion of Theorem 1.3 follows from Theorem 1.1 and the first part of
Theorem 1.3. \hfill\hb

\medskip

{\bf Proof of Theorem 1.5.}

\medskip

Here the arguments are similar to those in the proof of Theorem 1.1. So we only give those
proofs which is some what different and omitted other details.

First we assume that there exist only finitely many prime closed geodesics $\{c_k\}_{k=1}^q$
on a bumpy Finsler sphere $(S^d,F)$ with an odd integer $d\ge 2$ and assuming $i(c_k)\ge 2$
for $k=1,\ldots,q$.

\medskip

{\bf Claim 5.} {\it There exist at least $(d-1)$ distinct non-hyperbolic closed geodesics,
all of which possess even Morse indices, provided the total number of prime closed geodesics
on $(M,F)$ is finite.}

\medskip

In fact, firstly we can obtain the similar equations (\ref{4.12})-(\ref{4.14}) where $2N-1$
and $2N+1$ should be replaced by $2N-2$ and $2N+2$, respectively, i.e.,
for $1\le k\le q$, we obtain
\bea
&& i(c_k^{2m_k-m})\le 2N-2,\qquad \forall\ 1\le m\le 2m_k, \lb{4.40}\\
&& i(c_k^{2m_k})=2N-C(M_k)+2\Delta_k,\lb{4.41}\\
&& i(c_k^{2m_k+m})\ge 2N+2,\qquad \forall\ m\ge 1.\lb{4.42}\eea

Therefore, similarly to the equation (\ref{4.18}), we have
\be \sum_{p=0}^{2N+1}(-1)^p M_p = \sum_{k=1}^{q} 2m_k\ga_{c_k}
  -\sum_{1\le k\le q \atop i(c_k^{2m_k})\ge 2N+2} (-1)^{i(c_k^{2m_k})}
          \dim\ol{C}_{i(c_k^{2m_k})}(E,c_k^{2m_k}). \lb{4.43}\ee

Denote by $H_+^e,H_+^o,H_-^e,H_-^o$ the numbers similarly defined by (\ref{4.19})-(\ref{4.22}),
where $2N-1$ and $2N+1$ are replaced by $2N-2$ and $2N+2$ respectively. Note that now we can suppose
$N$ is a multiple of $D=d-1$.

Then by Claim 3, (\ref{4.43}), the definitions of $H^{e}_+$ and $H_+^o$, and Lemma 2.4, we have
\bea \frac{N(d+1)}{d-1}+H_+^o-H_+^e
&=& 2NB(d,1)+H_+^o-H_+^e \,=\, \sum_{k=1}^q 2m_k\gamma_{c_k}+H_+^o-H_+^e  \nn\\
&=& \sum_{p=0}^{2N+1}(-1)^p M_p \nn\\
&=& -(M_{2N+1}-M_{2N}+\cdots-M_1+M_0)\nn\\
&\le& -(b_{2N+1}-b_{2N}+\cdots-b_1+b_0)\nn\\
&=& \sum_{p=0}^{2N} b_p(\Lm M)\nn\\
&=& \frac{N(d+1)}{d-1}-\frac{d-1}{2},\lb{4.44}\eea
which yields
\be  H_+^e\ge H_+^e-H_+^o \ge \frac{d-1}{2}.\lb{4.45}\ee

Similarly, denote by $H_+^{'e},H_+^{'o},H_-^{'e},H_-^{'o}$ the numbers similarly defined
by (\ref{4.30})-(\ref{4.33}), where $2N'-1$ and $2N'+1$ are replaced by $2N'-2$ and $2N'+2$
respectively, and these numbers satisfy the following relationship
\be H_{\pm}^e = H_{\mp}^{'e},\qquad H_{\pm}^o = H_{\mp}^{'o}. \lb{4.46}\ee

Similarly to the inequality (\ref{4.36}), by the same arguments from (\ref{4.46}) we obtain
\be H_-^e = H_+^{'e} \ge H_+^{'e}-H_+^{'o} \ge \frac{d-1}{2}. \lb{4.47}\ee

Therefore from (\ref{4.45}) and (\ref{4.47}) we have
\be  q \ge H_+^e+H_-^e \ge \frac{d-1}{2}+\frac{d-1}{2} = d-1. \lb{4.48}\ee

Now by the same arguments in the proof of Theorem 1.1, it follows from the definitions of
$H_+^e$ and $H_-^e$ that these $(d-1)$ distinct closed geodesics are non-hyperbolic, and
the Morse indices of them and their iterations are all even. This proves Claim 5.

We denote these $(d-1)$ non-hyperbolic prime closed geodesics by $\{c_k\}_{k=1}^{d-1}$.

\medskip

{\bf Claim 6.} {\it There exist at least two distinct closed geodesics different from those
found in Claim 4 with even Morse indices provided the number of prime closed geodesics is finite.}

\medskip

In fact, for those $(d-1)$ distinct closed geodesics $\{c_k\}_{k=1}^{d-1}$ found in Claim 5,
there holds $i(c_k^{2m_k})\neq 2N$ by  the definitions of $H_+^e$ and $H_-^e$, which together
with (\ref{4.40}) and (\ref{4.42}) yields
\be  i(c_k^m)\neq 2N,\qquad m\ge 1,\quad k=1,\cdots,d-1.\lb{4.49}\ee
Then by Lemma 2.1 it yields
\be \sum_{1\le k\le d-1 \atop m\ge 1}\dim\ol{C}_{2N}(E,c_k^{m}) = 0. \lb{4.50}\ee

Therefore, noting that $N$ is a multiple of $d-1$, by (\ref{4.50}), the Morse inequality
(\ref{2.11}) of Theorem 2.5 and (\ref{2.5}) of Lemma 2.4, we obtain
\bea \sum_{d\le k\le q}\dim \ol{C}_{2N}(E,c_k^{2m_k})
&=& \sum_{d\le k\le q,\ m\ge 1}\dim \ol{C}_{2N}(E,c_k^m) \nn\\
&=& \sum_{1\le k\le q,\ m\ge 1}\dim \ol{C}_{2N}(E,c_k^m) \,=\, M_{2N} \,\ge\, b_{2N}(\Lm M) \,=\, 2,  \lb{4.51}\eea
where the first equality follows from (\ref{4.40}) and (\ref{4.42}).

Now by (\ref{4.51}) and Lemma 2.1, it yields that there exist at least two prime closed
geodesic $c_d$ and $c_{d+1}$ with $i(c_k^{2m_k})=2N$ and $i(c_k^{2m_k})-i(c_k)\in 2\N_0$ for
$k=d$ and $d+1$. Thus both $c_d$ and $c_{d+1}$ are different from $\{c_k\}_{k=1}^{d-1}$ by (\ref{4.49}),
and they and all of their iterates have even Morse indices. This completes the proof of Claim 6.

\medskip

Now Theorem 1.5 follows from Claim 5 and Claim 6. \hfill\hb

\medskip

{\bf Proof of Theorem 1.6.}

\medskip

In fact, under the curvature condition $\frac{\lm^2}{(1+\lm)^2}<K\le 1$, there holds
$i(c)\ge d-1\ge 2$ by Theorem 1 of \cite{Rad4} and Lemma 3 of \cite{Rad3}. Then Theorem 1.6 follows
from Theorem 1.5. \hfill\hb

\medskip

{\bf Proof of Theorem 1.8.}

\medskip

To prove Theorem 1.8, we assume that there exist only finitely many prime closed geodesics
$\{c_k\}_{k=1}^q$ on a compact simply connected bumpy Finsler manifold $(M,F)$ and all of
them are hyperbolic. Then we have $C(M_k)=0$ for all $1\le k\le q$. Note that (\ref{4.2})
yields $i(c_k^m)=mi(c_k)$ for all $m\ge 1$. So the positivity (\ref{4.3}) of the mean index
implies $i(c_k)\ge 1$ for all $1\le k\le q$.

Then by Theorem 3.5 and $i(c_k^{m+1})\ge i(c_k^m)$ for any $m\ge 1$, there exist infinitely many
$(q+1)$-tuples $(N, m_1, \cdots, m_q)\in\N^{q+1}$ such that for any $1\le k\le q$, there holds
\bea
i(c_k^{2m_k-m}) &=& 2N-i(c_k^m) \le 2N-1,\quad 1\le m<2m_k, \lb{4.52}\\
i(c_k^{2m_k}) &=& 2N,  \lb{4.53}\\
i(c_k^{2m_k+m}) &=& 2N+i(c_k^m) \ge 2N+1,\quad m\ge 1.\lb{4.54}\eea

By Lemma 2.1 and (\ref{4.54}), we know that all $c_k^{2m_k+m}$'s with $m\ge1$ and $1\le k\le q$ have no
contributions to the alternative sum $\sum_{p=0}^{2N}(-1)^p M_p$. Similarly also by Lemma 2.1 and
(\ref{4.52})-(\ref{4.53}), $i(c_k^{2m_k})$ and all $c_k^{2m_k-m}$'s with $1\le m< 2m_k$ and $1\le k\le q$
only have contribution to $\sum_{p=0}^{2N}(-1)^p M_p$.

Thus by Claim 3, Theorem 2.3, (\ref{4.23}) and (\ref{4.25}), we get
$$  -\frac{Ndn(n+1)}{D} = \sum_{k=1}^q 2m_k\gamma_{c_k}
= \sum_{p=0}^{2N}(-1)^p M_p \ge \sum_{p=0}^{2N}(-1)^p b_p = -\frac{Ndn(n+1)}{D}+ \frac{dn(n+1)}{4}, $$
which yields a contradiction. \hfill\hb

\medskip

{\bf Acknowledgements.} The authors would like to thank sincerely Dr. Hui Liu for his careful reading
of earlier versions of the manuscript. The first author also sincerely thank Professor Anatole Katok
for the support and hospitality during his visit to the Department of Mathematics of Pennsylvania State
University in August 2012 to August 2013.

\bibliographystyle{abbrv}

\begin{thebibliography}{100}\setlength{\itemsep}{-0.5mm}
\bibitem[Abr]{Abr} R. Abraham, {\it Bumpy metrics} in Global Analysis (Berkeley, 1968),
Proc. Sympos. Pure Math. 14, Amer. Math. Soc., Providence, 1968, 1-3.
\bibitem[Ano]{Ano} D. V. Anosov,  Gedesics in Finsler geometry. Proc. I.C.M.
(Vancouver, B.C. 1974), Vol. 2. 293-297 Montreal (1975) (Russian),
 {\it Amer. Math. Soc. Transl.} 109 (1977) 81-85.
\bibitem[BTZ1]{BTZ1} W. Ballmann, G. Thobergsson and W. Ziller, Closed geodesics
on positively curved manifolds. {\it Ann. of Math.} 116 (1982), 213-247.
\bibitem[BTZ2]{BTZ2} W. Ballmann, G. Thobergsson and W. Ziller, Existence of closed
geodesics on positively curved manifolds. {\it J. Diff. Geom.} 18 (1983), 221-252.
\bibitem[Ban1]{Ban1} V. Bangert, Geod\"{a}tische Linien auf Riemannschen Mannigfaltigkeiten,
{\it Jahresber. Deutsch. Math.-Verein.} 87 (1985), 39-66.
\bibitem[Ban2]{Ban2} V. Bangert, On the existence of closed geodesics on two-spheres.
{\it Inter. J. Math.} 4 (1993), 1-10.
\bibitem[BaK]{BaK} V. Bangert, W. Klingenberg, Homology generated by iterated closed geodesics.
{\it Topology} 22 (1983), 379-388.
\bibitem[BaL]{BaL} V. Bangert and Y. Long,   The existence of two closed
geodesics on every Finsler 2-sphere. {\it Math. Ann.}  346 (2010) 335-366.
\bibitem[BCS]{BCS} D. Bao, S. S. Chern and Z. Shen, An Introduction to
Riemann-Finsler Geometry. Springer. Berlin. 2000.
\bibitem[Bot]{Bot} R. Bott, On the iteration of closed geodesics and the Sturm
intersection theory. {\it Comm. Pure Appl. Math.}  9 (1956), 171-206.
\bibitem[Cha]{Cha} K. C. Chang, Infinite Dimensional Morse Theory and
Multiple Solution Problems. Birkh\"auser. Boston. 1993.
\bibitem[DuL1]{DuL1} H. Duan and Y. Long, Multiple closed geodesics on bumpy Finsler
$n$-spheres.  {\it J. Diff. Equa.} 233 (2007) 221-240.
\bibitem[DuL2]{DuL2} H. Duan and Y. Long, Multiplicity and stability of closed
geodesics on bumpy Finsler $3$-spheres. {\it Cal. Variations and PDEs}. 31 (2008)
483-496.
\bibitem[DuL3]{DuL3} H. Duan and Y. Long, The index growth and multiplicity
of closed geodesics. {\it J. Funct. Anal.} 259 (2010) 1850-1913.
\bibitem[DLW]{DLW} H. Duan, Y. Long and W. Wang, Two closed geodesics on compact
simply-connected bumpy Finsler manifolds. {\it J. Diff. Geom.}, to appear.
\bibitem[DLX]{DLX} H. Duan, Y. Long and Y. Xiao, Two closed geodesics on $\R P^{2n+1}$ with a
bumpy Finsler metric. {\it Calc. Variat. and PDEs.} 54 (2015), 2883-2894.
\bibitem[Fra]{Fra} J. Franks,  Geodesics on $S^2$ and periodic points of
annulus diffeomorphisms. {\it Invent. Math.}  108 (1992), 403-418.
\bibitem[GrM]{GrM} D. Gromoll and W. Meyer, Periodic geodesics on compact Riemannian
manifolds. {\it J. Diff. Geom.}  3 (1969), 493-510.
\bibitem[GuK]{GuK} J. Gutt, J. Kang, On the minimal number of periodic orbits on some
hypersurfaces in $\R^{2n}$. arXiv:1508.00166v1. (2015).
\bibitem[Hin]{Hin} N. Hingston,  Equivariant Morse theory and closed
geodesics. {\it J. Diff. Geom.} 19 (1984), 85-116.
\bibitem[HiR]{HiR} N. Hingston and H.-B. Rademacher, Resonance for loop homology of spheres.
{\it J. Diff. Geom.} 93 (2013) 133-174.
\bibitem[HWZ1]{HWZ1} H. Hofer, K. Wysocki and E. Zehnder, The dynamics on three-dimensional
strictly convex energy surfaces. {\it Ann. of Math.} 148 (1998) 197-289.
\bibitem[HWZ2]{HWZ2} H. Hofer, K. Wysocki and E. Zehnder,  Finite energy foliations
of tight three-spheres and Hamiltonian dynamics. {\it Ann. of Math.}  157 (2003), 125-257.
\bibitem[Kat]{Kat} A. B. Katok,  Ergodic properties of degenerate integrable
Hamiltonian systems. {\it Izv. Akad. Nauk SSSR.} 37 (1973) (Russian), {\it Math. USSR-Isv.}
 7 (1973), 535-571.
\bibitem[Kli]{Kli} W. Klingenberg, Lectures on Closed Geodesics. Springer. Berlin. 1978.
\bibitem[Liu]{Liu} C. Liu,  The relation of the Morse index of closed geodesics with
the Maslov-type index of symplectic paths. {\it Acta Math. Sinica.} 21 (2005), 237-248.
\bibitem[LLo]{LLo} C. Liu and Y. Long, Iterated Morse index formulae for closed
geodesics with applications. {\it Science in China}. 45. (2002) 9-28.
\bibitem[Lon1]{Lon1} Y. Long,  Bott formula of the Maslov-type index theory.
{\it Pacific J. Math.} 187 (1999), 113-149.
\bibitem[Lon2]{Lon2} Y. Long,  Precise iteration formulae of the Maslov-type index
theory and ellipticity of closed characteristics.  {\it Adv. Math.} 154 (2000), 76-131.
\bibitem[Lon3]{Lon3} Y. Long,  Index Theory for Symplectic Paths with Applications.
Progress in Math. 207, Birkh\"auser. 2002.
\bibitem[Lon4]{Lon4} Y. Long, Multiplicity and stability of closed geodesics on Finsler
2-spheres. {\it J. Euro. Math. Soc.} 8 (2006), 341-353.
\bibitem[LoD]{LoD} Y. Long and H. Duan,  Multiple closed geodesics on 3-spheres.
{\it Adv. Math.} 221 (2009) 1757-1803.
\bibitem[LoW]{LoW} Y. Long and W. Wang,  Stability of closed geodesics on Finsler
2-spheres. {\it J. Funct. Anal.} 255 (2008) 620-641.
\bibitem[LoZ]{LoZ} Y. Long and C. Zhu, Closed characteristics on compact convex hypersurfaces
in $\R^{2n}$. {\it Ann. of Math.} 155 (2002), 317-368.
\bibitem[Mor]{Mor1} M. Morse, Calculus of Variations in the Large. Amer.
Math. Soc. Colloq. Publ. vol. 18. Providence, R. I., Amer. Math. Soc. 1934.
\bibitem[Rad1]{Rad1} H.-B. Rademacher, On the average indices of closed geodesics. {\it J.
Diff. Geom.} 29 (1989), 65-83.
\bibitem[Rad2]{Rad2} H.-B. Rademacher, Morse Theorie und geschlossene Geodatische.
{\it Bonner Math. Schr.} 229 (1992).
\bibitem[Rad3]{Rad3} H.-B. Rademacher, A sphere theorem for non-reversible Finsler metric.
{\it Math. Ann.} 328 (2004), 373-387.
\bibitem[Rad4]{Rad4} H.-B. Rademacher, Existence of closed geodesics on positively curved
Finsler manifolds.  {\it Ergod. Th. \& Dynam. Sys.} 27 (2007), 957-969.
\bibitem[Rad5]{Rad5} H.-B. Rademacher, The second closed geodesic on the complex projective
plane. {\it Front. Math. China.}  3 (2008), 253-258.
\bibitem[Rad6]{Rad6} H.-B. Rademacher, The second closed geodesic on Finsler
spheres of dimension $n>2$. {\it Trans. Amer. Math. Soc. } 362 (2010), 1413-1421.
\bibitem[She]{She} Z. Shen, Lectures on Finsler Geometry. World Scientific.
Singapore. 2001.
\bibitem[ViS]{ViS1} M. Vigu\'e-Poirrier and D. Sullivan,  The homology theory of the closed
geodesic problem. {\it J. Diff. Geom.} 11 (1976), 633-644.
\bibitem[Wan1]{Wan1} W. Wang, Closed geodesics on positively curved Finsler spheres. {\it Adv.
Math.} 218 (2008) 1566-1603.
\bibitem[Wan2]{Wan2} W. Wang, On a conjecture of Anosov. {\it Adv. Math.} 230 (2012), 1597-1617.
\bibitem[Wan3]{Wan3} W. Wang, Non-hyperbolic closed geodesics on Finsler spheres. {\it J. Diff. Geom.}
99 (2015), 473-496.
\bibitem[XiL]{XiL} Y. Xiao and Y. Long, Topological structure of non-contractible loop space
and closed geodesics on real projective spaces with odd dimensions. {\it Adv. Math.} 279
(2015), 159-200.
\bibitem[Zil]{Zil} W. Ziller,  Geometry of the Katok examples. {\it Ergod. Th. \& Dynam. Sys.
3 (1982), 135-157. }


\end{thebibliography}

\end{document}